\documentclass{gSTA2e}

\usepackage{epstopdf}
\usepackage{subfigure}

\usepackage{amsbsy}
\usepackage{amsfonts}
\usepackage{amsmath}
\usepackage{amssymb}
\usepackage{amsthm}
\usepackage{units}
\usepackage[inline]{enumitem}
\usepackage{textcomp}

\usepackage{epsfig}
\usepackage{graphicx}
\usepackage{natbib}
\usepackage{rotating}

\usepackage[T1]{fontenc}
\usepackage[english,brazilian]{babel}
\usepackage[utf8]{inputenx}

\usepackage{setspace}
\usepackage[hyphenbreaks]{breakurl}
\usepackage{pdflscape}
\usepackage{icomma}
\usepackage{glossaries-preamble} 
\usepackage[capitalise]{cleveref} 
\DeclareMathOperator{\sign}{sign}

\theoremstyle{plain}
\newtheorem{theorem}{Theorem}[section]
\newtheorem{coro}[theorem]{Corollary}

\newtheorem{prop}[theorem]{Proposition}

\crefformat{prop}{#2Prop.~#1#3}
\Crefformat{prop}{#2Prop.~#1#3}
\crefformat{coro}{#2Cor.~#1#3}
\Crefformat{coro}{#2Cor.~#1#3}
\crefmultiformat{coro}{Cors.~#1#1#3}{ and~#2#1#3}{, #2#1#3}{ and~#2#1#3}
\Crefmultiformat{coro}{Cors.~#2#1#3}{ and~#2#1#3}{, #2#1#3}{ and~#2#1#3}
\crefformat{section}{#2Sect.~#1#3}
\Crefformat{section}{#2Sect.~#1#3}
\crefformat{table}{#2Table~#1#3}
\Crefformat{table}{#2Table~#1#3}
\crefformat{figure}{#2Fig.~#1#3}
\Crefformat{figure}{#2Fig.~#1#3}
\crefmultiformat{figure}{Figs.~#2#1#3}{ and~#2#1#3}{, #2#1#3}{ and~#2#1#3}
\Crefmultiformat{figure}{Figs.~#2#1#3}{ and~#2#1#3}{, #2#1#3}{ and~#2#1#3}
\crefrangeformat{figure}{Figs.~#3#1#4--#5#2#6}
\Crefrangeformat{figure}{Figs.~#3#1#4--#5#2#6}
\crefformat{equation}{#2Eq.~(#1)#3}
\Crefformat{equation}{#2Eq.~(#1)#3}
\crefmultiformat{equation}{Eqs.~(#2#1#3)}{ and~(#2#1#3)}{, (#2#1#3)}{ and~(#2#1#3)}
\Crefmultiformat{equation}{Eqs.~(#2#1#3)}{ and~(#2#1#3)}{, (#2#1#3)}{ and~(#2#1#3)}
\crefrangeformat{equation}{Eqs.~(#3#1#4)--(#5#2#6)}
\Crefrangeformat{equation}{Eqs.~(#3#1#4)--(#5#2#6)}

\theoremstyle{definition}

\theoremstyle{remark}

\begin{document}


\articletype{Original Paper}

\title{Closed-form mathematical expressions for the \al{ECR} distribution}

\selectlanguage{brazilian}
\author{
  \name{Thiago VedoVatto\textsuperscript{a}$^{\ast}$\thanks{$^\ast$Corresponding author. Email: thiago.vedovatto@ifg.edu.br}
  and Abraão David Costa do Nascimento\textsuperscript{b}}
  \affil{\textsuperscript{a}Instituto Federal de Educação, Ciência e Tecnologia de Goiás, Goiânia, Brasil;
  \textsuperscript{b}Universidade Federal de Pernambuco, Recife, Brasil}
  \received{\today}
}
\selectlanguage{english}

\maketitle

\begin{abstract}
  The \af{CR} distribution has been successfully used to describe asymmetric and heavy-tail events from radar imagery.
  Employing such model to describe lifetime data may then seem attractive, but some drawbacks arise: its \al{pdf} does not cover non-modal behavior as well as the \as{CR} \af{hrf} assumes only one form.
  To outperform this difficulty, we introduce an extended \as{CR} model, called \af{ECR} distribution.
  This model has two parameters and \as{hrf} with decreasing, decreasing-increasing-decreasing and upside-down bathtub forms.
  In this paper, several closed-form mathematical expressions for the \as{ECR} model are proposed: median, mode, probability weighted, log-, incomplete and order statistic moments and \al{fei}.
  We propose three estimation procedures for the \as{ECR} parameters: \af{ML}, bias corrected \as{ML} and \al{PB} methods.
  A simulation study is done to assess the performance of estimators.
  An application to survival time of heart problem patients illustrates the usefulness of the \as{ECR} model.
  Results point out that the \as{ECR} distribution may outperform classical lifetime models, such as the gamma, \al{BSd}, Weibull and \al{LN} laws, before heavy-tail data.
\end{abstract}


\begin{classcode}62E99; 62P10; 62N02; 62F10\end{classcode}

  \section{Introduction}


  \Ab{SAR} systems have been indicated as important tools to solve remote sensing issues.
  Among other, this fact may be justified by their capability in operating on all weather conditions and in providing high resolution images.
  However, \ab{SAR} images are strongly contaminated by an interference pattern called \emph{speckle noise}.
  Thus, working with \ab{SAR} imagery requires a specialized modeling.
  The \ab{CR} distribution, known as generalized \al{Cauchy} \citep{BibalanAmindavar-NonGaussianamplitude-2016}, has received great attention as descriptor of \ab{SAR} features.
  In this paper, we propose a new distribution that extends the \ab{CR} law as describing for lifetime data.

  \citet{KuruogluZerubia-ModelingSARimages-2004} presented evidence that this law may outperform the Weibull, \ab{LN} and \ab{KBessel} distributions.
  \citet{PengZhao-SARImageFiltering-2014} used a mixture of \abp{CR} as an approximation of the \ab{HTR} distribution in order to model \ab{SAR} features.
  \citet{LiEkman-CauchyPowerAzimuth-2010,LiEkman-CauchyRayleighscattering-2011} introduced the \ab{CR} law as a model of scattering clusters in state space based on simulation model for single and multiple polarization channels.
  \citet{HillAchimBullAl-Mualla-Dualtreecomplex-2014} produced a novel bivariate shrinkage technique to provide a quantitative improvement in image denoising using the \ab{CR} model.
  Recently, \citet{BibalanAmindavar-newalphaand-2015,BibalanAmindavar-NonGaussianamplitude-2016} furnished a mathematical treatment to the \ab{HTR} distribution by means of mixtures of \ab{CR} and \ab{Rayleigh} models.
  This approach was developed for modeling amplitude of ultrasound images through the \ab{HTR} distribution.


  In this paper, we have twofold goal.
  First using the original and a new two-parameter extended \ab{CR} distributions in the survival analysis context as alternatives to the gamma, \ab{BSd}, Weibull, \ab{LN} models and other lifetime ones.
  Several works have used mixtures to model heavy-tail lifetimes \citep{HarrisonMillard-Balancingacuteand-1991,McCleanMillard-Modellinginpatient-1993,GardinerLuoTangRamamoorthi-FittingHeavyTailed-2014}.
  But this modeling scheme becomes the associated estimation process hard.
  The reason of applying the \ab{CR} law for describing lifetimes is firstly to accommodate heavy-tail distributions.
  However, its \ab{pdf} does not cover non-modal behavior (as the exponential does) as well as the \ab{CR} \ab{hrf} assumes only one form, limiting its use in practice.
  To that end, we extend the \ab{CR} distribution  by using the exponentiation.
  \citet{Gompertz-NatureFunctionExpressive-1825} and \citet{Verhulst-Noticesurla-1838} pioneered this method.
  It may be described briefly as: Let $X$ be a \ab{rv} with \ab{cdf} $G(\cdot)$, its exponentiated version is a \ab{rv} with \ab{cdf} $G(\cdot)^\beta$, where $\beta>0$ is an additional shape parameter.
  This technique has triggered several models: some examples are the \ab{EE} \citep{GuptaKundu-Generalizedexponentialdistributions-1999}, exponentiated Rayleigh (two-parameter \ab{BX}) \citep{KunduRaqab-GeneralizedRayleighdistribution-2005,SurlesPadgett-Somepropertiesscaled-2005}, \ab{EHC} \citep{CordeiroLemonte-BetaHalfCauchy-2011} and \ab{ELi} \citep{NadarajahBakouchTahmasbi-generalizedLindleydistribution-2011} laws.
  We denote the new model as \ab{ECR} distribution.
  This model has two parameters and accommodates \abp{hrf} with decreasing, decreasing-increasing-decreasing and upside-down bathtub forms.
  As second aim, several closed-form mathematical expressions for the \ab{ECR} model are proposed: median, mode, probability weighted, log-, incomplete and order statistic moments and \ab{fei}.
  There are several concepts for the term ``closed-form''; as one example, one has the work of \citet{StoverWeisstein-ClosedFormSolution-2017}.
  Here, we understand ``closed-form'' as
  \begin{quotation}
    ``all expressions in terms of a finite number of known (special) functions, which have well-defined analytic properties.''
  \end{quotation}
  Some of these quantities do not exist for specific \ab{ECR} parametric regions, as it will be discussed in our presentation.
  In additional, we propose two estimation procedures for the \ab{ECR} parameters: \abp{MLE} and \abp{PBE}.
  It is known \abp{MLE} present a bias of order $O(n^{-1})$ for small and moderate sample sizes  $(n)$, where $O(\cdot)$ is the Landau notation to represent order.
  To overcome it we furnish a second-order expression for the \ab{ECR} bias according to \citet{CoxSnell-GeneralDefinitionResiduals-1968} and, consequently, propose a second bias-corrected \ab{MLE}.
  A Monte Carlo simulation study is made to quantify the performance of proposed estimators, adopting bias and \ab{SSD} as figures of merit.
  An application to survival time of patients in a waiting list of heart transplant is performed to illustrate the usefulness of the \ab{ECR} model.
  Results point out that the \ab{ECR} distribution may outperform well-defined biparametric models; such as-beyond other seventeen ones-the gamma, \ab{BSd}, Weibull and \ab{LN} laws.


  This paper is outlined as follows.
  In \cref{sec:Thealpha-stablemodel}, we present a brief discussion about the \ab{CR} model (special case of the \ab{HTR} law) obtained from the $\alpha$-stable distribution.
  The proposed model is presented in \cref{sec:The broached model} and some of its descriptive properties are listed in \cref{sec:ECR:SomeElementaryProperties}.
  \Cref{sec:ECR:Moments} discusses some \ab{ECR} mathematical expressions.
  \Cref{sec:ECR:Pointestimation} addresses procedures to obtain the \abp{MLE}, \abp{PBE} and \abp{CS-MLE}.
  A simulation study is presented in \cref{sec:ECR:SimulationStudy}.
  A real data application is yielded in \cref{sec:ECR:Application}.
  Finally, \cref{sec:ECR:Concludingremarks} provides concluding remarks.

  \section{The \al{aS}, \al{HTR} and \al{CR} models}\label{sec:Thealpha-stablemodel}

  \citet{Levy-Calculdesprobabilites-1925} pioneered the \ab{aS} distribution, which has been widely used in financial time series \citep{Mandelbrot-ParetoLevyLaw-1960,FamaSchwert-Humancapitaland-1977,Voit-TheStatisticalMechanicsofFinancialMarkets-2005} and, recently, applied in \ab{SAR} image processing \citep{Pierce-RCScharacterizationusing-1996,WangLiaoLi-ShipDetectionin-2008,PengXuZhouZhao-SARimageclassification-2011}.
  This law has not tractable pdf expression and is often represented by its \ab{cf} given by: For $-\infty<t<\infty$, 
  \begin{equation*}
    \sym{cfS}_{\alpha S}(t)=
    \begin{cases}
      \exp\left\{j\zeta t-\lambda|t|^\alpha\left[1+j \rho \sign(t)\tan\left( \frac{\alpha\pi}{2}\right)\right]\right\},\quad\text{if }\alpha\ne 1,\\
      \exp\left\{j\zeta t-\lambda|t|\left[1+j \rho \sign(t) \frac{2}{\pi} \log|t|\right]\right\},\quad\text{if }\alpha= 1,
    \end{cases}
  \end{equation*}
  where $-\infty<\zeta<\infty$, $\lambda>0$, $0<\alpha\le 2$ and $-1\le\rho\le 1$ are the location, scale, characteristic and symmetry parameters, respectively.

  Now consider the (zero-mean) \ab{SaS} distribution having \ab{cf}
  \begin{equation}\label{eq:S-alpha-Scharacteristicfunction}
    \sym{cfS}_{S\alpha S}(t) = \exp(-\lambda|t|^\alpha).
  \end{equation}
  According to \citet{KuruogluZerubia-ModelingSARimages-2004}, one can define the bivariate \ab{IaS} distribution, which has \ab{cf} as
  \begin{equation}\label{eq:S-alpha-Sisotropic}
    \sym{cfS}_{I\alpha S}(t_{\sym{Re}},t_{\sym{Im}})=\exp(-\lambda|\bm{t}|^\alpha),
  \end{equation}
  where $t_{\sym{Re}}$ and $t_{\sym{Im}}$ can be understood as outcomes of real and imaginary parts of a complex \ab{rv}, say $ \bm{t} $, respectively, and $|\bm{t}|=\sqrt{t_{\sym{Re}}^2+t_{\sym{Im}}^2}$ represents the amplitude of $\bm{t}$.
  From \cref{eq:S-alpha-Scharacteristicfunction}, the bivariate \ab{IaS} \ab{pdf} can be obtained by taking the Fourier transform of \cref{eq:S-alpha-Sisotropic}:
  \begin{equation*}
    f_{\alpha,\lambda}(x_{\sym{Re}},x_{\sym{Im}})= \frac{1}{(2\pi)^2} \int_{t_{\sym{Re}}}\int_{t_{\sym{Im}}}\!\exp(-\lambda|\bm{t}|^\alpha)\exp[ -j2\pi(x_{\sym{Re}}t_{\sym{Re}}+x_{\sym{Im}}t_{\sym{Im}}) ]\,\mathrm{d}t_{\sym{Im}}\,\mathrm{d}t_{\sym{Re}}.
  \end{equation*}
  Representing this integral into its polar form in terms of $s=|\bm{t}|$ and $\omega=\arctan( \nicefrac{t_{\sym{Re}}}{t_{\sym{Im}}})$, one has
  \begin{equation}\label{eq:polarform}
    f_{\alpha,\lambda}(x_{\sym{Re}},x_{\sym{Im}})= \frac{1}{(2\pi)^2} \int_{0}^{2\pi} \int_{0}^{\infty}\!s\exp(-\lambda s^\alpha)J_0(s|\bm{x}|)\,\mathrm{d}s\,\mathrm{d}\omega,
  \end{equation}
  where $J_0$ is the zero order Bessel function of the first kind \citep{AbramowitzStegun-Handbookofmathematicalfunctions-1972} and $|\bm{x}|=\sqrt{x^2_{\sym{Re}}+x^2_{\sym{Im}}}$.
  Since $\omega$ does not appear in the integrand of the \cref{eq:polarform}, it collapses in
  \begin{equation}\label{eq:polarformintegrated}
    f_{\alpha,\lambda}(x_{\sym{Re}},x_{\sym{Im}})= \frac{1}{2\pi} \int_{0}^{\infty}\!s\exp(-\lambda s^\alpha)J_0(s|\bm{x}|)\,\mathrm{d}s.
  \end{equation}
  Now consider to determine the \ab{pdf} of the corresponding amplitude. 
  Using polar coordinates transformation (where $r$ and $\phi$ indicate amplitude and phase of $x_{\sym{Re}}+ \bm{j} x_{\sym{Im}}$ such that $\bm{j}=\sqrt{-1}$, respectively), the following joint \ab{pdf} can be obtained: For $r>0$ and $0\le\phi\le2\pi$,
  \begin{equation}\label{eq:polarcoordinates}
    f(r,\phi)=rf_{\alpha,\lambda}(r\cos(\phi),r\sin(\phi)).
  \end{equation}

  Under the common assumption that $\phi$ is uniformly distributed on $[0,2\pi]$ \citep{KuruogluZerubia-ModelingSARimages-2004}, we can determine an expression for the amplitude distribution from replacing \eqref{eq:polarformintegrated} in \eqref{eq:polarcoordinates} and integrating over $\phi$: 
  \begin{equation}\label{eq:amplitudemodel}
    g_{\alpha,\lambda}(r)=r \int_{0}^{\infty}\!s\exp(-\lambda s^\alpha) J_0(rs)\,\mathrm{d}s.
  \end{equation}
  This expression is know as the \ab{HTR} \ab{pdf} \citep{NikiasShao-SignalProcessingwithAlpha-1995}.

  When $\alpha=1$ in \eqref{eq:amplitudemodel} and using the identity 6.623.2 given by \citet{GradshteynRyzhik-Tableofintegralsseriesandproducts-2007}, the \ab{CR} \ab{pdf} is obtained: For $r>0$,
  \begin{equation*}
    g_{\as{CR}}(r) = \frac{\lambda r}{(\lambda^2+r^2)^{\nicefrac{3}{2}}},
  \end{equation*}
  where $\lambda>0$ is the scale parameter.
  This \ab{pdf} is associated with the amplitude of a coefficient on which the components are jointly \al{Cauchy} distributed \citep{KuruogluZerubia-ModelingSARimages-2004}.
  The \ab{CR} \ab{cdf} is given by
  \begin{equation}\label{eq:CauchyRayleighcdf}
    G_{\as{CR}}(r) = 1- \frac{\lambda}{\sqrt{\lambda^2+r^2}}.
  \end{equation}


  Under a non-physical perspective, some remarks are important.
  The \ab{CR} distribution is a case nested in the \ab{CW} distribution \citep[p.~157]{Rinne-TheWeibullDistributionAhandbook-2009} with \ab{cdf} (for $x>0$)
  \begin{equation}\label{eq:compoundWeibull}
    F_{\ab{CW}}(x)=1-\exp(-ax^c)\left[1+\left(\frac{x}{\lambda}\right)^c\right]^{-k },
  \end{equation}
  where $ c\ge0 $ and $k\ge0$ are shape parameters and $a\ge0$ and $\lambda>0$ are scale parameters. 
  The \ab{CR} distribution is obtained at $ c=2 $, $ k=\nicefrac{1}{2} $ and $a=0$.
  At $a=0$, the \ab{CW} distribution coincides with the \ab{SM}, \ab{aka} three parameter \ab{BXII}, model, which is another notable \ab{CR} extension.  
  The \ab{CR} distribution is also a simpler case of other models, like the \ab{FP} \citep{Feller-AnIntroductiontoProbabilityandItsApplications-2-1971,ArnoldLaguna-OnGeneralizedParetoDistributionswithApplicationstoIncomeData-1977} and \ab{TB} \citep{KlugmanPanjerWillmot-LossModelsFromDatatoDecisions-WileySeriesinProbabilityandStatistics-2012} laws.
  However, although this uniparametric model has been successfully used as reviewed above, it does not seem to exist works about mathematical properties of possible \ab{CR} biparametric extensions.
  In what follows, we cover this gap, presenting several new closed-form properties of the \ab{ECR} model.
  Some \ab{ECR} asymptotic theory results are also provided.

  \section{The broached model}\label{sec:The broached model}

  In this section, we aim to introduce a \ab{CR} extension as alternative to classical biparametric models (as the gamma and Weibull distributions), which we denote as \ab{ECR} model.
  The \ab{ECR} distribution has \ab{cdf} (for $x>0$) given by
  \begin{equation}\label{eq:ECR:cdf}
    F(x)=\left(1- \frac{\lambda}{\sqrt{\lambda^2+x^2}}\right)^\beta,
  \end{equation}
  where $\lambda>0$ is a scale parameter and $\beta>0$ is an additional shape parameter and, as a consequence, its \ab{pdf} (for $x>0$) is
  \begin{equation}\label{eq:ECR:pdf}
    f(x)=\beta\lambda \frac{x}{(\lambda^2+x^2)^{\nicefrac{3}{2}}}\left(1- \frac{\lambda}{\sqrt{\lambda^2+x^2}}\right)^{\beta-1}.
  \end{equation}
  From \cref{eq:ECR:cdf,eq:ECR:pdf}, the \ab{ECR} \ab{hrf} is expressed as
  \begin{equation}
    \sym{hrfS}(x) = \beta  \lambda  \frac{x}{\left(\lambda ^2+x^2\right)^{\nicefrac{3}{2}}} \left(1-\frac{\lambda }{\sqrt{\lambda ^2+x^2}}\right)^{\beta -1} \left[1-\left(1-\frac{\lambda }{\sqrt{\lambda ^2+x^2}}\right)^{\beta }\right]^{-1}\label{eq:ECR:hrf}.
  \end{equation}
  We denote this case as $ X\sim\operatorname{\as{ECR}}(\beta,\lambda) $.  

  The \ab{ECR} law is a model nested in the \ab{EBXII} \citep{Al-HussainiAhsanullah-ExponentiatedDistributions-AtlantisStudiesinProbabilityandStatistics5-2015}, \ab{KLL} \citep{SantanaOrtegaCordeiroSilva-KumaraswamyLogLogistic-2012}, \ab{GFP} \citep{Zandonatti-DistribuzionideParetoGeneralizzate-UniversityofTrento-2001}, \ab{KBXII} \citep{ParanaibaOrtegaCordeiroPascoa-KumaraswamyBurrXII-2013} and \ab{McBXII} \citep{Mead-NewGeneralizationBurr-2014,GomesSilvaCordeiro-TwoExtendedBurr-2015} distributions.
  In general, the properties obtained for these models do not have closed-form expressions, which may imply in computational difficulties and an insufficient mathematical treatment (where divergence cases may happen).
  As previously discussed, the \ab{ECR} model can be very useful, but some of its properties do not has closed-form expressions and, therefore, a specialized treatment is required.
  We do it in the next sections.

  The \ab{ECR} \ab{qf} is obtained by inverting \eqref{eq:ECR:cdf}: For $0<p<1$,
  \begin{equation}\label{eq:ECR:quantilefunction}
    \sym{qfS}(p)= \frac{\lambda}{1-p^{\nicefrac{1}{\beta}}}\sqrt{2p^{\nicefrac{1}{\beta}}-p^{\nicefrac{2}{\beta}}}.
  \end{equation}
  The \ab{ECR} median is then determined from using $p=\nicefrac{1}{2}$ in \eqref{eq:ECR:quantilefunction}:
  \begin{equation*}
    \sym{median}=\frac{\lambda}{2^{\nicefrac{1}{\beta}}-1}\sqrt{2^{\frac{\beta+1}{\beta}}-1}.
  \end{equation*}
  We can also use \eqref{eq:ECR:quantilefunction} for generating outcomes of a \ab{ECR} \ab{rv} by the \ab{itsm} as follows:
  \begin{itemize}
    \item Generate $u$ as an outcome of $U\sim U(0,1)$;
    \item Obtain $x=\sym{qfS}(u)$ as an outcome from the \ab{ECR} distribution.
  \end{itemize}

\begin{figure}
    \centering
    \subfigure[\as{ECR} \as{pdf} and sample histogram midpoints]{
      \includegraphics[width=0.45\textwidth,trim=0 15 0 55,clip]{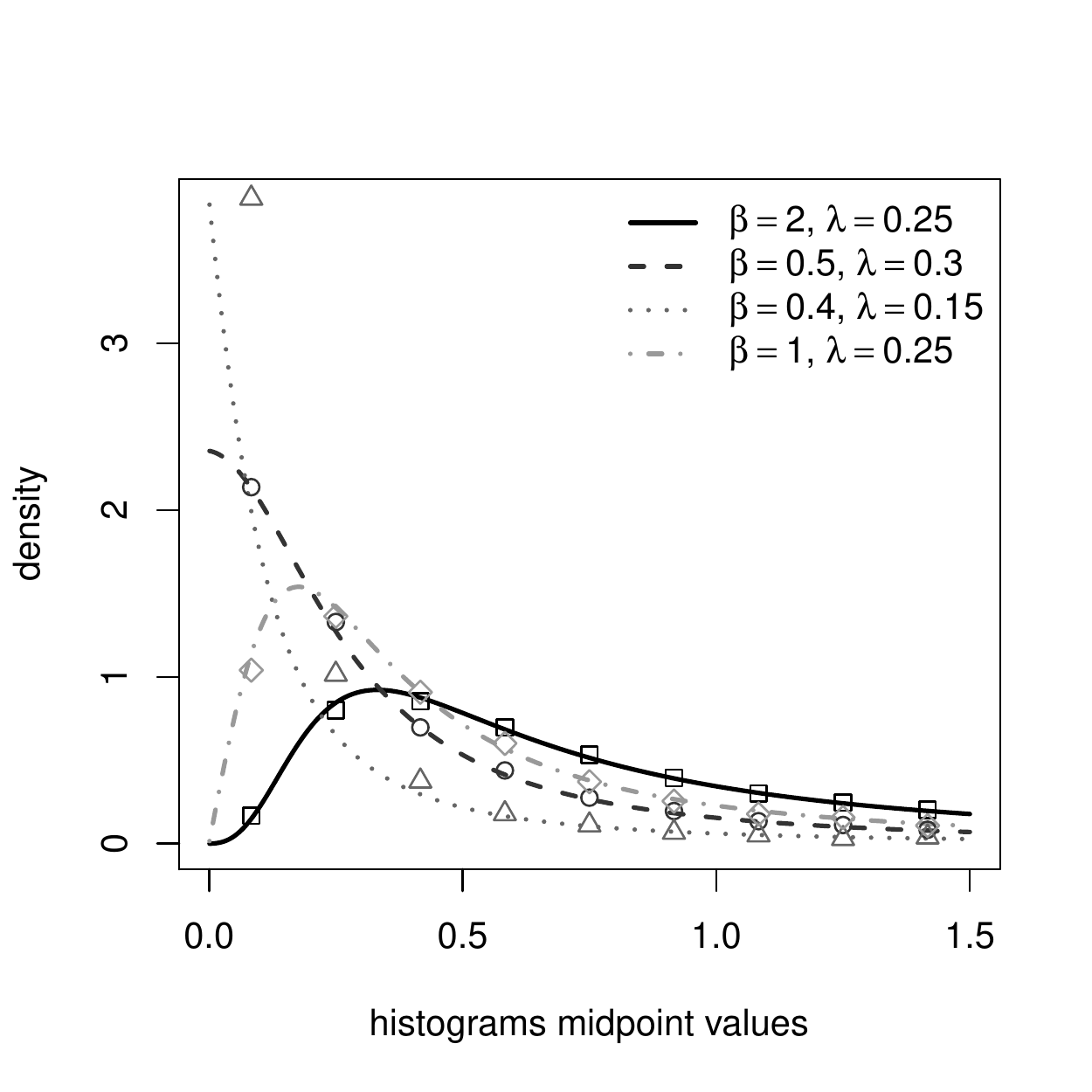}
      \label{fig:ECR:DensityPlot}
    }
    \subfigure[\as{ECR} \as{hrf}]{
      \includegraphics[width=0.45\textwidth,trim=0 15 0 55,clip]{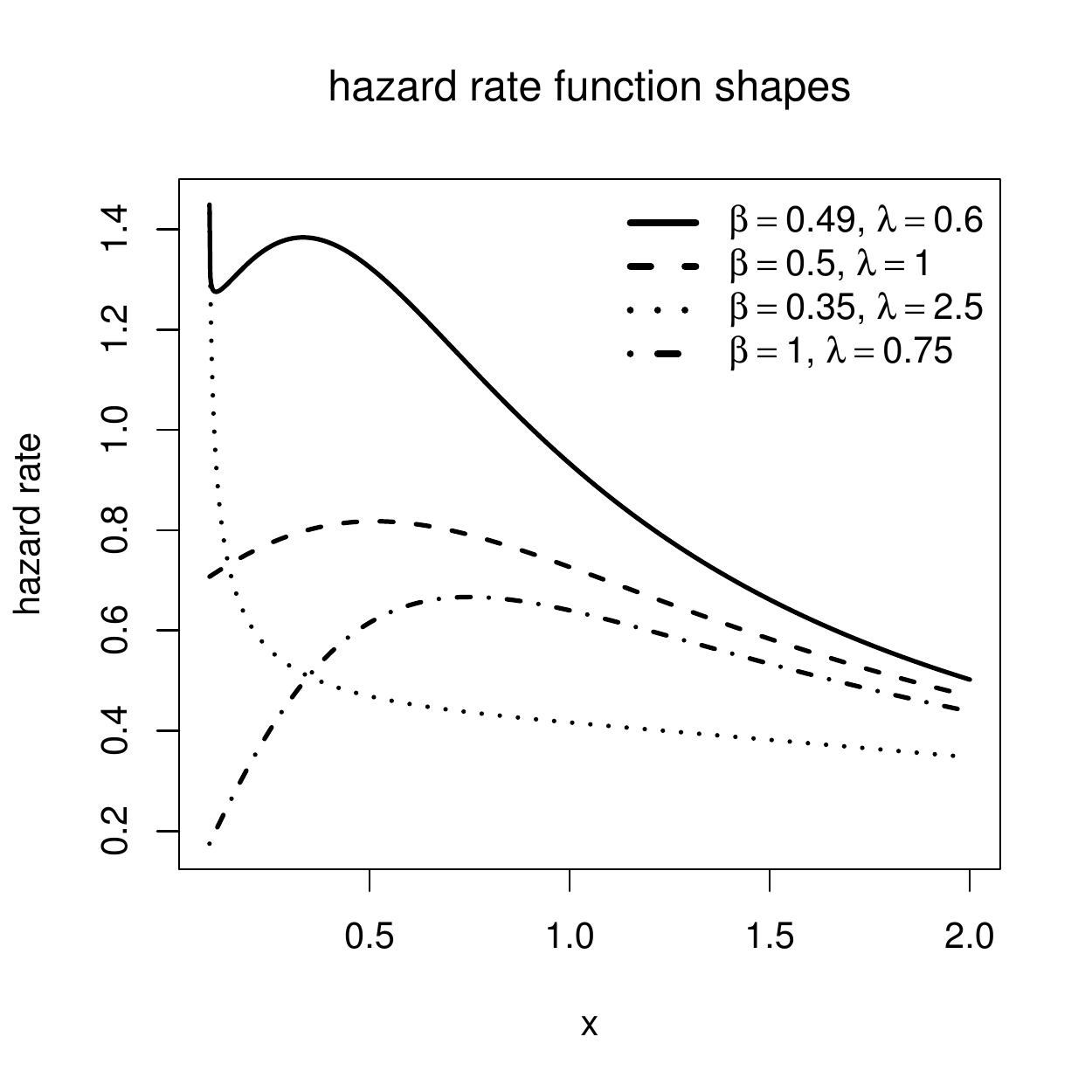}
      \label{fig:ECR:HazardRatePlot}
    }
    \caption{\as{ECR} \as{hrf} and \as{pdf} plots for various values of $\beta$ and $\lambda$.}
    \label{fig:ECR:pdfandhrf}
  \end{figure}

  \Cref{fig:ECR:DensityPlot,fig:ECR:HazardRatePlot} depict some possible shapes of \eqref{eq:ECR:pdf} and \eqref{eq:ECR:hrf} for some parameter values.
  The \ab{ECR} \ab{pdf} may be unimodal and non-modal decreasing from both infinite or a finite value, property detailed subsequently.
  It is also noticeable that the \ab{ECR} \ab{hrf} has three basic shapes: decreasing, decreasing-increasing-decreasing and upside-down bathtub. 
   In contrast, \ab{CR} \ab{hrf} only assumes upside-down bathtub shape (starting at the origin). 

  \section{Some elementary properties}\label{sec:ECR:SomeElementaryProperties}

  In this section we derive some \ab{ECR} descriptive properties.
  \Cref{prop:ECR:limitspdf} shows that the \ab{ECR} \ab{pdf} asymptote is more flexible than the \ab{CR} one.
  The \ab{ECR} \ab{pdf} limit assumes three distinct values when $x$ approaches zero (in terms of $\beta$), while it is always null for the \ab{CR} case.

  \begin{prop}\label{prop:ECR:limitspdf}
    The \ab{ECR} \ab{pdf} has the following limiting behavior:
    \begin{equation}\label{eq:ECR:limitspdf}
      \lim_{x\to 0^+}f(x)=
      \begin{cases}
        \infty,                   & \beta<\nicefrac{1}{2}, \\
        \frac{\sqrt{2}}{2\lambda}, & \beta=\nicefrac{1}{2}, \\
        0,                        & \beta>\nicefrac{1}{2}.
      \end{cases}
    \end{equation}
  \end{prop}

  A \ab{rv} $ X $ is called \ab{rvi} with tail index $ a>0 $ if $ \sym{sfS}_X(cx)\sim c^{-a}\sym{sfS}_X(x) $ as $ x\to\infty $. 
  \citet[p.~55]{CookeNieboerMisiewicz-RegularlyVaryingand-2014}, citing \citet{Chistyakov-TheoremSumsIndependent-1964}, defined the \emph{subexponential} distributions containing the \ab{rvi} class.
  \citet[Lemma~3.2]{FossKorshunovZachary-AnIntroductiontoHeavy-TailedandSubexponentialDistributions-2011}, crediting \citet{Chistyakov-TheoremSumsIndependent-1964} again, pointed out that any subexponential distribution with positive support is \emph{long-tailed} and therefore it is \emph{heavy-tailed} too.
  Thus the \ab{ECR} model belongs to all these classes.
  \Cref{prop:ECR:rvi} assigns the \ab{ECR} distribution to the \ab{rvi} class.

  \begin{prop}\label{prop:ECR:rvi}
    The \ab{ECR} distribution is \ab{rvi} with tail index $ a=1 $.
  \end{prop}

  \Cref{prop:ECR:ECRmode} shows a closed-form expression for the \ab{ECR} mode. 
  The \ab{ECR} \ab{pdf} is non-modal for $\beta<\nicefrac{1}{2}$ as induced from \cref{prop:ECR:limitspdf} and non-modal for $\beta=\nicefrac{1}{2}$.
  \begin{prop}\label{prop:ECR:ECRmode}
    The \ab{ECR} mode is given by
    \begin{equation}\label{eq:ECR:ECRmode}
      \sym{mode}=\frac{\lambda}{2\sqrt{2}}\sqrt{(\beta+1)^2+(\beta-1)\sqrt{\beta^2+6\beta+17}},\quad\beta>\nicefrac{1}{2}.
    \end{equation}
  \end{prop}

  \begin{coro}\label{coro:ECR:ECRmodellimits}
    The \ab{ECR} mode has the following limiting behavior $\lim_{\beta\to\nicefrac{1}{2}^+}\sym{mode}=0$.
  \end{coro}

  \begin{coro}\label{coro:ECR:CRmodel}
    The \ab{CR} mode is expressed by $ \sym{mode}_{\as{CR}}=\nicefrac{\lambda}{\sqrt{2}} $.
  \end{coro}

\begin{figure}[htb!]
    \centering
    \subfigure[mode]{
      \includegraphics[width=0.30\textwidth,trim=40 60 45 65,clip]{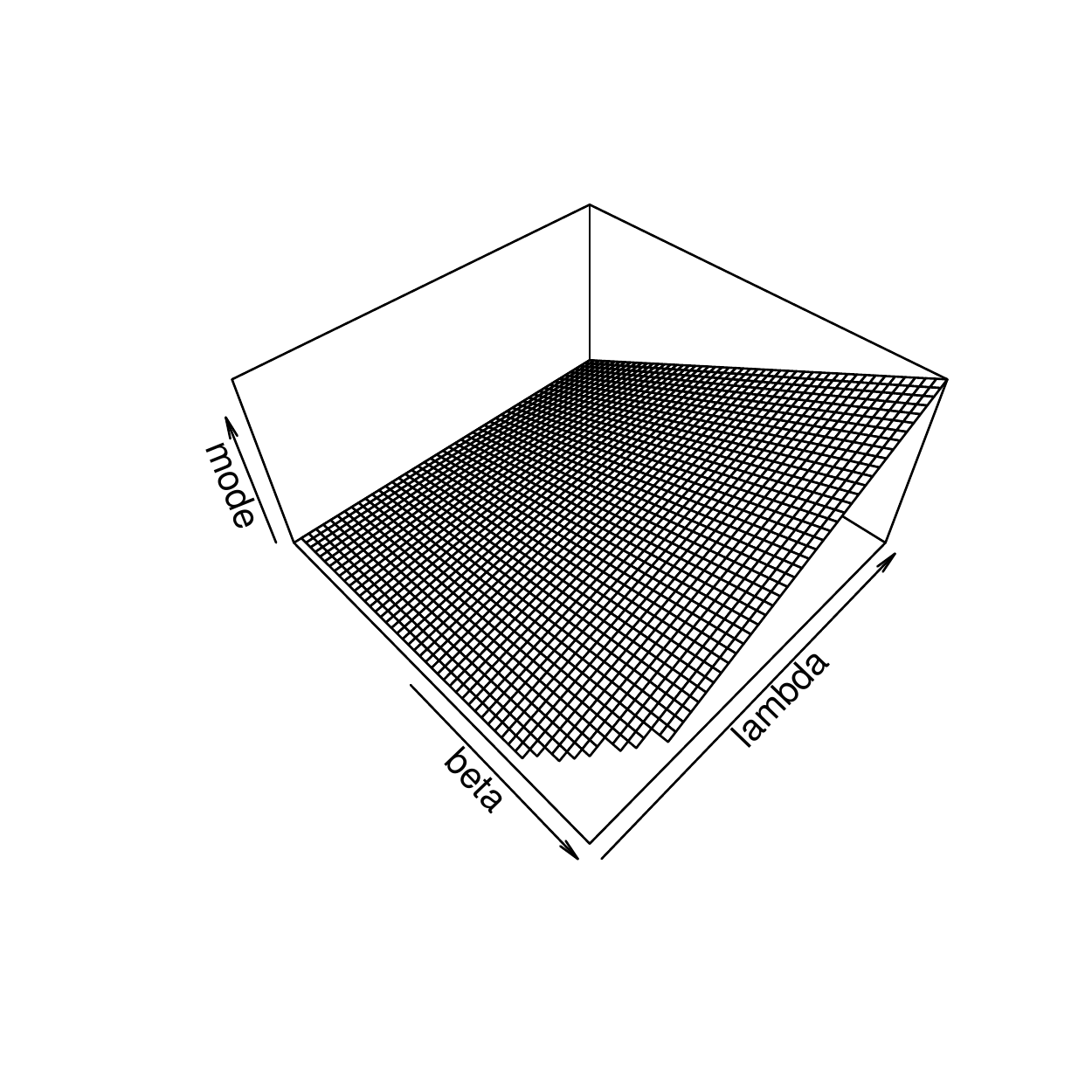}
      \label{fig:ECR:modesurface}
    }
    \subfigure[\ab{median}]{
      \includegraphics[width=0.30\textwidth,trim=40 60 45 65,clip]{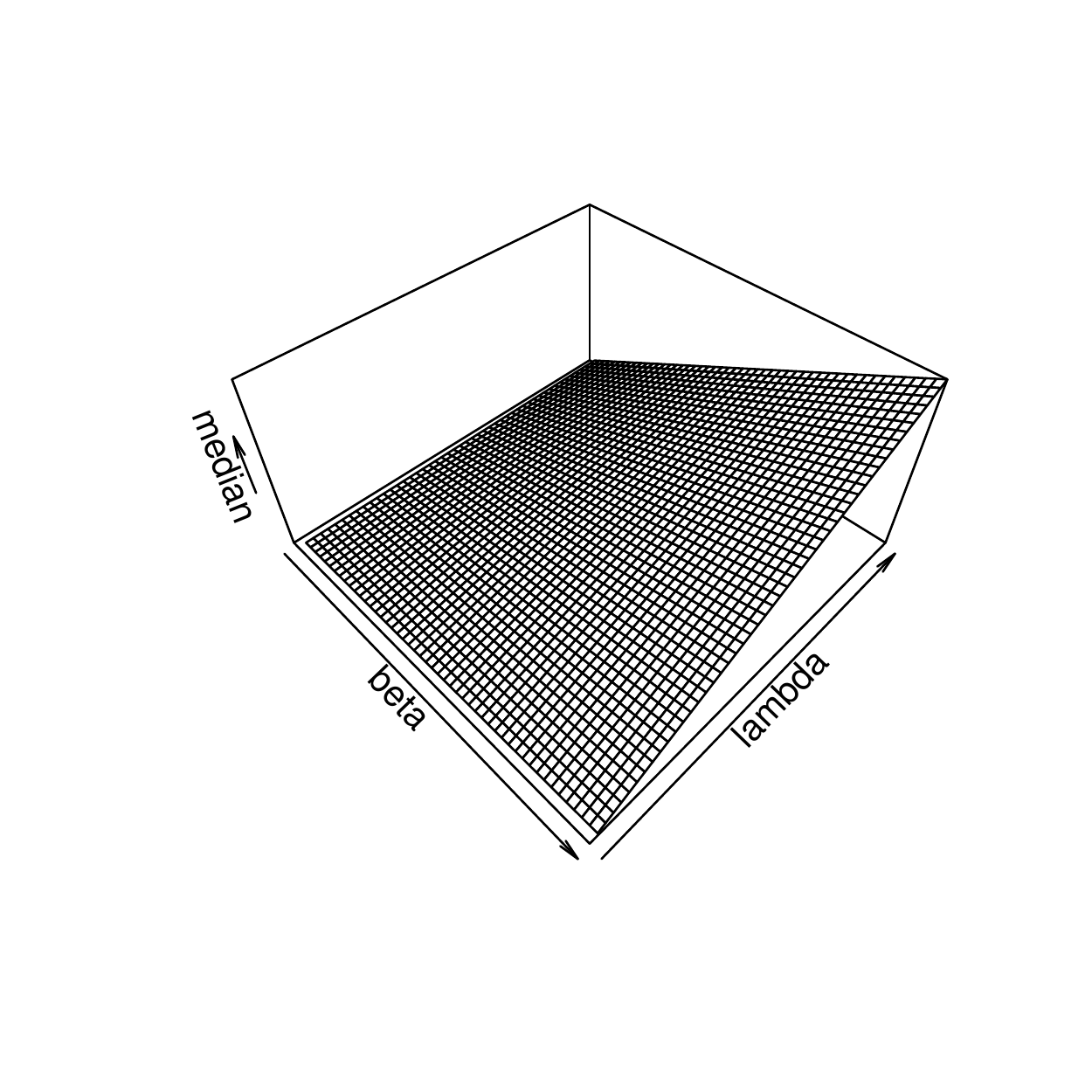}
      \label{fig:ECR:mediansurface}
    }
    \subfigure[\ab{median} - mode]{
      \includegraphics[width=0.30\textwidth,trim=40 60 45 65,clip]{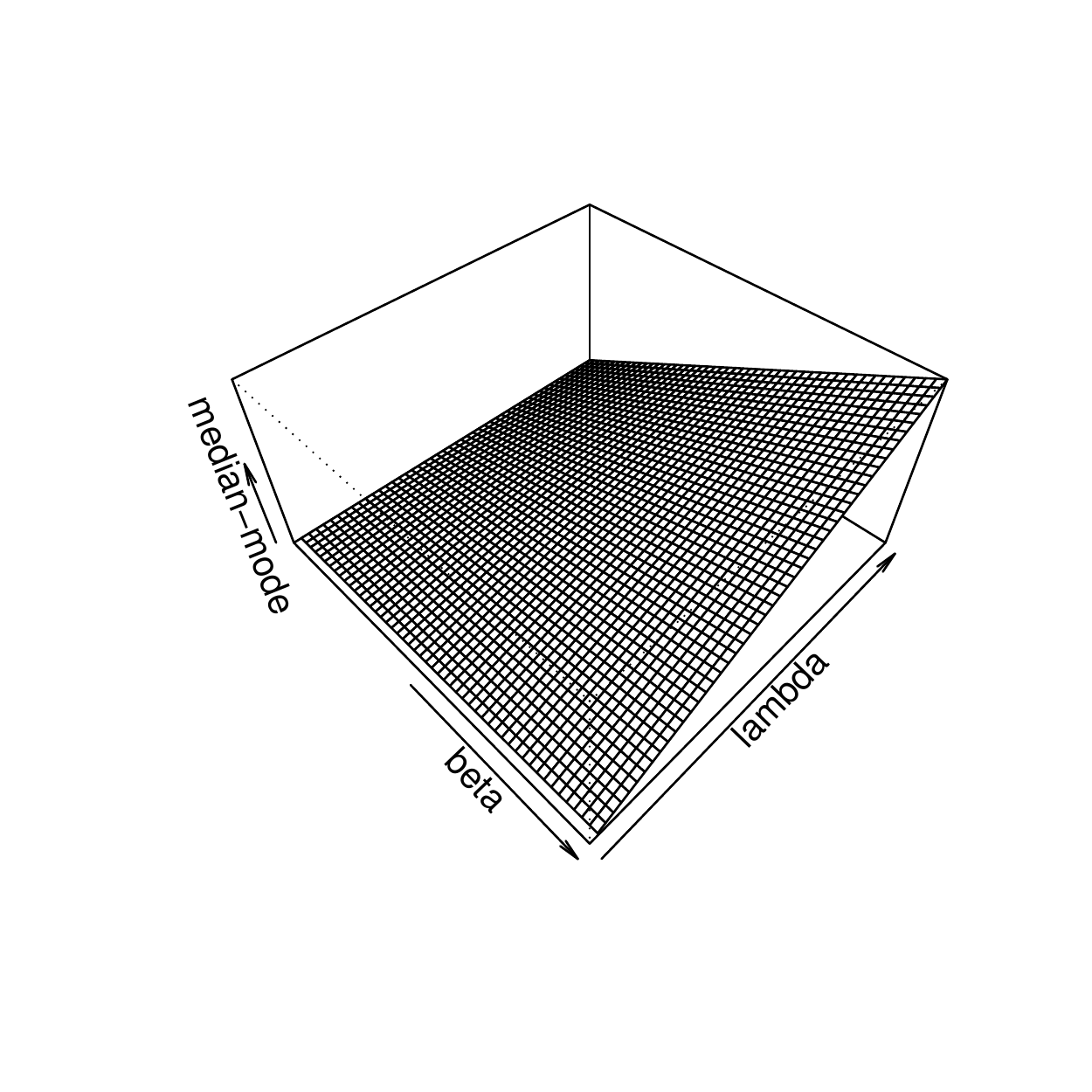}
      \label{fig:ECR:diffsurface}
    }
    \caption{Plots of surface of \as{ECR} mode and median and its differences.}
    \label{fig:ECR:modemediancomparison}
  \end{figure}

  \Cref{fig:ECR:modemediancomparison} exhibits plots for \ab{ECR} mode and \ab{median} and $\sym{median}-\sym{mode}$.
  \cref{fig:ECR:diffsurface} indicates that our model has always positive skewness.

  \Cref{prop:ECR:auxprop1} shows the \ab{ECR} distribution belongs to the scale family.
  Two immediate consequences of this proposition are $\sym{E}(X^k)=\lambda^k\sym{E}(Z^k) \label{eq:ECR:EXk}$ and $\sym{E} (\log X) =\log\lambda+\sym{E} (\log Z)$.
  In the rest of this paper, we denote $Z\sim\operatorname{\as{ECR}}(\beta,1)$ as the standard \ab{ECR} \ab{rv}.
  \begin{prop}\label{prop:ECR:auxprop1}
    If $Z$ is the standard \ab{ECR} \ab{rv}, then $X=\lambda Z\sim\operatorname{\as{ECR}}(\beta,\lambda)$.
  \end{prop}

  \section{Moments}\label{sec:ECR:Moments}

  Expressions for the $ r $-moments, \abp{im} and \abp{osm} of some \ab{ECR} extensions were explored by \citet{ParanaibaOrtegaCordeiroPescim-betaBurrXII-2011,ParanaibaOrtegaCordeiroPascoa-KumaraswamyBurrXII-2013}, \citet{Al-HussainiAhsanullah-ExponentiatedDistributions-AtlantisStudiesinProbabilityandStatistics5-2015}, \citet{SilvaGomes-SilvaRamosCordeiro-ExponentiatedBurrXII-2015}, \citet{Mead-NewGeneralizationBurr-2014,GomesSilvaCordeiro-TwoExtendedBurr-2015} and \citet{CordeiroOrtegaHamedaniGarcia-McBurrXIIand-2016}.
  But all results in these papers were not given in closed-form.
  In general, they were deduced from power series expansions.
  In this section, we show that several moment \ab{ECR} expressions often do not exist and, therefore, it is required a detailed discussion about it.

  \subsection{Probability weighted moments}\label{sec:ECR:ProbabilityWeightedMoments}

  The \abp{pwm} with order indexes $ r, s, t\in\sym{Re} $, \sym{pwmS}, for a \ab{rv} $X\sim\operatorname{\as{ECR}}(\beta,\lambda)$ are defined by \glsname{pwmS}.
  \Cref{prop:ECR:pwms} presents a closed-form expression for the \ab{ECR} \abp{pwm} restricting $ t\in\sym{Z} $.

  \begin{prop} \label{prop:ECR:pwms}
    Let $X\sim\operatorname{\as{ECR}}(\beta,\lambda)$, for $-2(s+1)\beta<r<1$ and $ t\in\sym{Z} $,
    \begin{align*}\label{eq:pwmsECR}
      \sym{pwmS} = & \beta(\lambda\sqrt{2})^r\,
      \sum_{i=0}^{t}(-1)^t\binom{t}{i}\,
      \sym{bf}\left(1-r,\frac{r}{2}+(s+i+1)\beta \right)\,\nonumber\\
      &\times\sym{2F1}\left(-\frac{r}{2},\frac{r}{2}+(s+i+1)\beta;1-\frac{r}{2}+(s+i+1)\beta;\,\frac{1}{2}\right),
    \end{align*}
    where $\sym{2F1}$ is the \al{2F1} defined as
    \begin{equation*}
      \glsname{2F1},\quad \glsname{Pochhammer},	
    \end{equation*}
    $\sym{bf}(\cdot,\cdot)$ is the \ab{bf} and $ \sym{gf}(\cdot) $ is the \ab{gf}.
  \end{prop}

  \citet[p.~33]{FossKorshunovZachary-AnIntroductiontoHeavy-TailedandSubexponentialDistributions-2011} showed that in the \ab{rvi} class with index $ a>0 $ all moments of order $ 0<r<a $ are finite, while all moments with order $ r>a $ are infinite.
  This fact occurs in the \ab{ECR} distribution, and they will be presented in \cref{coro:ECR:momentsECR}. Closed-form expressions for the \ab{ECR} moments consist in a simple product involving the beta and the \abp{2F1}.

  \begin{coro}\label{coro:ECR:momentsECR}
    Let $X\sim\operatorname{\as{ECR}}(\beta,\lambda)$, for $r\in(-2\beta,1)$,
    \begin{equation}\label{eq:ECR:momentsECR}
      \sym{E}(X^r) = \beta(\lambda\sqrt{2})^r\,
      \sym{bf}\left(1-r,\frac{r}{2}+\beta \right)\,
      \sym{2F1}\left(-\frac{r}{2},\frac{r}{2}+\beta;1-\frac{r}{2}+\beta;\,\frac{1}{2}\right).
    \end{equation}
  \end{coro}

  \Cref{coro:ECR:momentsCR} presents the \ab{CR} moments, derived from the \cref{coro:ECR:momentsECR}.
  The obtained result agrees with the known results for the \ab{SM} distribution as presented by \citet[eq.~6.46]{KleiberKotz-StatisticalSizeDistributionsinEconomicsandActuarialSciences-2003} (Assuming $ a=2 $, $b=\lambda$ and $ q=\nicefrac{1}{2} $ in the \ab{BXII} parameterization adopted by the referred authors).
  \begin{coro}\label{coro:ECR:momentsCR} 
    Let $Y\sim\operatorname{\as{CR}}(\lambda)$, for $ r\in(-2,1) $
    \begin{equation}\label{eq:momentsCR}
      \sym{E}(Y^r)= \frac{\lambda^r}{\sqrt{\pi}}\,
      \sym{gf}\left( \frac{1-r}{2}\right)\,
      \sym{gf}\left(1+ \frac{r}{2}\right)
      =\frac{\lambda^r}{2}\,
      \sym{bf}\left(\frac{1-r}{2},1+\frac{r}{2}\right).
    \end{equation}
  \end{coro}

  \subsection{First log-moment and an induced regression model}
  \Cref{prop:ECR:ElogX} contains an useful result to define a log-linear regression model or a location quantity (since the mean is not  finite) for the \ab{ECR} distribution.
  Some integrals used in the proof of \cref{prop:ECR:ElogX} will be useful for the next results.
  \begin{prop}\label{prop:ECR:ElogX}
    Let $ X\sim\operatorname{\as{ECR}}(\beta,\lambda) $, then
    \begin{equation*}
      \sym{E}(\log X) = \log\lambda + 
      \frac{1}{2}\sym{lpf}\left(\frac{1}{2};1,\beta\right) + 
      \sym{dgf}(1+\beta)+\sym{EM} - 
      \frac{1}{\beta},
    \end{equation*}
    where $ \sym{lpf}(\cdot;\cdot,\cdot) $ is the \ab{lpf}, $ \sym{dgf}(\cdot) $ is the \ab{dgf} and $ \sym{EM}=0.57721\ldots $ is the \ab{EM} defined by, respectively:
    \begin{equation*}
      \glsname{lpf},\quad\glsname{dgf}\text{ and }\glsname{EM}.
    \end{equation*}
  \end{prop}

  An outline of a possible \ab{ECR} regression is given as follows.
  Let $\epsilon_i \sim \operatorname{\as{ECR}}(\alpha,1)$ and ${\bm z}_i = (1, z_{i1}, \ldots, z_{ip})^ \top$ be a vector of predict variables.  
  Assuming that there is interest in describing linearly the expected value of logarithm of a \ab{rvi} positive variable, the following result holds:
  \begin{equation*}
    Y_i = \mathrm{e}^{{\bm z}_i^\top{\bm \beta}}\epsilon_i\sim \operatorname{\as{ECR}}(\alpha,\exp(\bm{z}_i^\top\bm{\beta})) \quad \Longleftrightarrow \quad
    \log(Y_i)  = {\bm z}_i^\top{\bm \beta} + \log(\epsilon_i),
  \end{equation*}
  where ${\bm \beta}_0^* \mathrel{\mathop:}= \beta_0 +\sym{E} [-\log(\epsilon_i)]$ and ${\bm \beta}=(\beta_0^*, \beta_1, \ldots, \beta_p)^\top$ is a vector of regression coefficients.
  In this case, notice that
  $\sym{E}(\log Y_i)=\beta_0+\sum_{k=1}^p\beta_kz_{ik}$
  and it is possible to proof that $\sym{Var}(\log Y_i)=\sym{Var}(\log\epsilon_i)<\infty$.
  We omit the expression of $\sym{Var}({\log Y_i})$ by simplicity, but it may be obtained from author's contact.

  \subsection{Incomplete moments}
  Let $X$ be an \ab{ECR} \ab{rv}, its \abp{im} are defined by \Cref{prop:ECR:ims}.
  It is known \abp{im} consist in the main part of important inequality tools in income applications \citep{ButlerMcDonald-Usingincompletemoments-1989} such as \citeauthor{Lorenz-MethodsMeasuringConcentration-1905} curves \citep{Lorenz-MethodsMeasuringConcentration-1905} and \citeauthor{Gini-MeasurementInequalityIncomes-1921} measure \citep{Gini-MeasurementInequalityIncomes-1921}.
  For future applications, users of the \ab{ECR} law may employ the \cref{prop:ECR:ims} to obtain these measures.
  \begin{prop}\label{prop:ECR:ims}
    Let $ X\sim\operatorname{\as{ECR}}(\beta,\lambda) $, for $ r>-2\beta $
    \begin{equation*}
      \sym{imS}(x_0)=\frac{\beta2^{\nicefrac{r}{2}+1}\lambda^ru_0^{\beta+\nicefrac{r}{2}}}{2\beta+r}\sym{F1}\left(\frac{r}{2}+\beta,r,-\frac{r}{2};\frac{r}{2}+\beta+1;u_0,\frac{u_0}{2}\right),
    \end{equation*}
    where $ u_0=1-\nicefrac{\lambda}{\sqrt{x_0^2+\lambda^2}} $ and \sym{F1} is the \ab{F1} defined by
    \begin{equation*}
      \glsname{F1}.
    \end{equation*}
  \end{prop}

  \subsection{Order statistics}
  Moments of the order statistics play an vital role in quality control and reliability issues \citep{AhsanullahNevzorovShakil-AnIntroductiontoOrderStatistics-2013}, where researchers wish to predict the future failure items based on recorded failure times.
  Let $ X_1,\ldots,X_n $ be a random sample from the \ab{ECR} law and \glsname{os} be the corresponding order statistics.
  Let \sym{ospdfS} be the \ab{ospdfS} \sym{os}:
  \begin{equation*}
    \glsname{ospdfS},
  \end{equation*}
  where $ \sym{cdfS}(x) $ and $ \sym{pdfS}(x) $ are given in \eqref{eq:ECR:cdf} and \eqref{eq:ECR:pdf}, respectively.
  Then we can express the \ab{ECR} order statistics \ab{pdf} as
  \begin{equation*}
    \sym{ospdfS}(x)=\frac{\beta\lambda}{\sym{bf}(i,n-i+1)} \frac{x}{(\lambda^2+x^2)^{\nicefrac{3}{2}}}\left(1-\frac{\lambda}{\sqrt{\lambda^2+x^2}}\right)^{i\beta-1}\left[1-\left(1-\frac{\lambda}{\sqrt{\lambda^2+x^2}}\right)^{\beta}\right]^{n-i}.
  \end{equation*}
  In \cref{prop:ECR:ros}, we derive a closed-form expression for the \ab{ECR} \abp{osm}.
  \begin{prop}\label{prop:ECR:ros}
    Let \sym{os} be the \ab{os} of a $n$-points random sample from $ X\sim\operatorname{\as{ECR}}(\beta,\lambda) $, for $ r\in(-2i\beta,1) $,
    \begin{align*}
      \sym{E}(X_{i:n}^r)=&\frac{\beta(\lambda\sqrt{2})^r}{\sym{bf}(i,n-i+1)} \sum_{j=0}^{n-i}(-1)^j\binom{n-i}{j}
      \sym{bf}\left(1-r,\frac{r}{2}+(i+j)\beta\right)\nonumber\\
      &\times\sym{2F1}\left(\frac{-r}{2},\frac{r}{2}+(i+j)\beta;1-\frac{r}{2}+(i+j)\beta;\frac{1}{2}\right).
    \end{align*}
  \end{prop}

  \section{Inference and second-order asymptotic theory under the \as{ECR} model}\label{sec:ECR:Pointestimation}

  \subsection{\Alp{MLE}}\label{sec:Maximum likelihood estimation}

  In this section, we provide a procedure to find \abp{MLE} for \ab{ECR} parameters as well as construct asymptotic confidence intervals and hypothesis tests.
  Let $x_1,\ldots,x_n$ be an observed sample obtained from the \ab{ECR} distribution with vector parameter $\sym{PV}=(\beta,\lambda)^\top$.
  The \ab{lf} at $\bm{\theta}$ is expressed by
  \begin{equation*}
    \sym{lfS}(\sym{PV})= \prod_{i=1}^{n}\beta\lambda \frac{x_i}{(\lambda^2+x_i^2)^{\nicefrac{3}{2}}}\left[1- \frac{\lambda}{\sqrt{\lambda^2+x_i^2}}\right]^{\beta-1}
  \end{equation*}
  and the corresponding \ab{llf} is
  \begin{equation}\label{eq:ECR:log-likelihood}
    \sym{llfS}(\sym{PV})    = n\log(\beta\lambda)
    +\hspace{-1ex}\underbrace{T_1(\lambda,\bm{x})}_{\displaystyle\sum_{i=1}^{n}\log x_i}
    \hspace{-1ex}+\hspace{0.5ex}3\hspace{-6.5ex}\overbrace{T_2(\lambda,\bm{x})}^{\displaystyle\sum_{i=1}^{n} \log\left(\frac{1}{\sqrt{\lambda^2+x_i^2}}\right)}\hspace{-6.5ex}
    +\hspace{0.5ex}(\beta-1)\hspace{-8.5ex}\underbrace{T_3(\lambda,\bm{x})}_{\displaystyle\sum_{i=1}^{n}\log\left(1- \frac{\lambda}{\sqrt{\lambda^2+x_i^2}}\right)}\hspace{-8.5ex}.
  \end{equation}
  The \ab{ML} estimate for $\sym{PV}$ is the pair that maximizes the \ab{llf}; i.e.,
  $
  \widehat{\sym{PV}}=\arg\max_{\sym{PV}\in\sym{PS}}[\sym{llfS}(\sym{PV})]
  $,
  where $\sym{PS}$ is the parametric space.
  \abp{MLE} have not closed-form and, therefore, it is required using iterative numerical methods such as Newton-Raphson and BFGS, implemented in several softwares and programming languages (\as{eg}, \texttt{R}, \texttt{Julia}, \texttt{Mathematica} and \texttt{Maple}). 

  Other way to find ML estimates is by \ab{llf} derivatives, known as \emph{score functions}.
  The components of the \ab{SV} $ \sym{SV}=\left(\sym{SC}_\beta,\sym{SC}_\lambda\right)^\top$ are
  \begin{equation}\label{eq:scorealpha}                                                      \sym{SC}_\beta = \frac{n}{\beta}+ T_3(\lambda,\bm{x}),
  \end{equation}
  and $ \sym{SC}_\lambda =\nicefrac{n}{\lambda}+(1-\beta)T_4(\lambda,\bm{x})-(\beta+2)T_5(\lambda,\bm{x}) $, where $T_4(\lambda,\bm{x})= \sum_{i=1}^{n} \nicefrac{1}{\sqrt{\lambda^2+x_i^2}}$ and $T_5(\lambda,\bm{x})= \lambda \sum_{i=1}^{n} \nicefrac{1}{(\lambda^2+x_i^2)}$.
  The \ab{ML} estimates, $\boldsymbol{\hat{\theta}}=(\hat{\beta},\hat{\lambda})^\top$, can also be obtained numerically by solving the nonlinear equations system $ \sym{SC}_\beta = \sym{SC}_\lambda = 0 $.

  From \cref{eq:scorealpha}, it is possible to obtain a semi-closed \ab{MLE} for $\beta$:
  \begin{equation*}
    \hat{\beta}(\lambda)=-\frac{n}{T_3(\lambda,\bm{x})}.
  \end{equation*}
  By replacing $\beta$ by $\hat{\beta}(\lambda)$ in \eqref{eq:ECR:log-likelihood}, the following \ab{pllf} for $\lambda$ is obtained:
  \begin{equation*}
    \sym{llfS}_{\beta}(\lambda)= n\left[\log\left( \frac{-n\lambda}{T_3(\lambda,\bm{x})}\right)-1\right] + T_1(\lambda,\bm{x})+ 3T_2(\lambda,\bm{x})-T_3(\lambda,\bm{x}).
  \end{equation*}
  We can also obtain the \ab{ML} estimates for $\lambda$ by maximizing the \ab{pllf} $\sym{llfS}_{\beta}(\lambda)$ with respect to $\lambda$.
  The \ab{psf} can be expressed as
  \begin{equation*}
    \sym{SC}_{\lambda|\beta=\hat{\beta}(\lambda)}= \frac{n}{\lambda}
    +\left(1+\frac{n}{T_3(\lambda,\bm{x})}\right)T_4(\lambda,\bm{x})
    -\left(2-\frac{n}{T_3(\lambda,\bm{x})}\right)T_5(\lambda,\bm{x}).
  \end{equation*} 
  The \ab{ML} estimate for $\lambda$, $\widehat{\beta}(\widehat{\lambda})$, can also be defined as a root of the equation $ \sym{SC}_{\lambda|\beta=\hat{\beta}(\lambda)}=0 $.  

  \subsection{The \al{fei}}
  Now consider to obtain the \ab{fei} for the \ab{ECR} model.
  Among its potentialities, \ab{fei} allows to determine the asymptotic confidence interval based on the following result \citep[p.~386]{BickelDoksum-MathematicalStatisticsBasicIdeasandSelectedTopics-1-2000}:
  \begin{equation*}
    \sqrt{n}(\hat{\sym{PV}}-\sym{PV})\sim \sym{N}_2(0,\sym{feiS}^{-1}(\sym{PV}))\text{ as }n\to\infty,
  \end{equation*}
  where \sym{feiS} is the \ab{fei} defined by
  \begin{equation}\label{eq:ECR:FisherInformationMatrix}
    \sym{feiS}=\sym{feiS}(\sym{PV})=-\frac{1}{n}
    \begin{bmatrix} 
      \sym{cum}_{\beta\beta} &  \sym{cum}_{\beta\lambda} \\
      \sym{cum}_{\lambda\beta} & \sym{cum}_{\lambda\lambda}
    \end{bmatrix}.
  \end{equation}
  The entries of $\sym{feiS}$ obey the notation:
  \begin{equation*}
    \sym{SCD}=\dfrac{\partial\sym{llfS}(\sym{PV})}{\partial\theta_1\cdots\partial\theta_n}\text{ and }\sym{cumnlld}=\sym{E}\left(\sym{SC}_{\theta_1\cdots\theta_n}\right).
  \end{equation*}

  The \ab{fei} is obtained in \cref{prop:ECR:FEI}. 
  In general, determining this matrix requires hard integrations, but we provide solutions for all in the \ab{ECR} case.
  \citet{MahmoudEl-Ghafour-FisherInformationMatrix-2015} obtained the \ab{fei} for the \ab{GFP} distribution and it is also possible to obtain the \ab{ECR} \ab{fei} from their expressions.

  \begin{prop}\label{prop:ECR:FEI}
    The \ab{ECR} \ab{fei} elements in \eqref{eq:ECR:FisherInformationMatrix} are given by
    \begin{subequations}
      \begin{align}
        \sym{cum}_{\beta\beta} & = -\frac{n}{\beta^2} \label{eq:ECR:FEIMele1}, \\
        \sym{cum}_{\beta\lambda} & = \sym{cum}_{\lambda\beta}= \frac{n}{\lambda}\left(\frac{2}{\beta +2 }-\frac{3}{\beta +1 }\right)\label{eq:ECR:FEIMele2} \intertext{and}
        \sym{cum}_{\lambda\lambda} & = \frac{n}{\lambda^2}\left(\frac{18}{\beta +2}-\frac{36}{\beta +3}+\frac{16}{\beta +4}-1\right)\label{eq:ECR:FEIMele3}.
      \end{align}
    \end{subequations}
  \end{prop}

  \subsection{\alp{CS-MLE}}
  \Cref{coro:ECR:invFEI,coro:ECR:FEIder} and \cref{prop:ECR:TD} are the basic results to construct the \abp{CS-MLE}. \Cref{coro:ECR:invFEI} presents the inverse \ab{fei}, which is also useful for determining analytic expressions for the asymptotic standard errors of the \abp{MLE} and construct confidence intervals.
  \begin{coro}\label{coro:ECR:invFEI}
    The inverse \ab{fei} is given by	
    $$
    \sym{feiS}^{-1}=
    \sym{feiS}^{-1}(\sym{PV})=-n
    \begin{bmatrix}
      \sym{cum}^{\beta,\beta}   & \sym{cum}^{\beta,\lambda} \\
      \sym{cum}^{\lambda,\beta} & \sym{cum}^{\lambda,\lambda}
    \end{bmatrix},
    $$
    where $ \sym{cum}^{\beta,\beta} = -\nicefrac{1}{n}\left[\frac{\beta ^2 (\beta +1)^2 (\beta +2) \left(\beta ^2+11 \beta +36\right)}{\beta ^3-7 \beta ^2+10 \beta +72}\right] $, $ \sym{cum}^{\beta,\lambda} =\sym{cum}^{\lambda,\beta}= \nicefrac{\lambda}{n}\left[\frac{\beta  (\beta +1) (\beta +2) (\beta +3) (\beta +4)^2 }{\beta ^3-7 \beta ^2+10 \beta +72}\right] $ and $ \sym{cum}^{\lambda,\lambda} = -\nicefrac{\lambda ^2}{n}\left[\frac{(\beta +1)^2 (\beta +2)^2 (\beta +3) (\beta +4)}{\beta  \left(\beta ^3-7 \beta ^2+10 \beta +72\right)}\right] $.
  \end{coro}

  \Cref{coro:ECR:FEIder} presents the first derivatives of the \ab{fei} components.
  We will assume the notation $\sym{dercum}=\nicefrac{\partial\sym{cum}_{\theta_i\theta_j}}{\partial\theta_k}$ for the derivatives of \ab{fei} components.
  \begin{coro}\label{coro:ECR:FEIder}
    The first derivatives of the \ab{fei} components are given by $\sym{cum}_{\beta\beta}^{(\beta)} =\nicefrac{2 n}{\beta ^3}$, $\sym{cum}_{\beta\beta}^{(\lambda)} =0$, $\sym{cum}_{\beta\lambda}^{(\beta)} =\nicefrac{n}{\lambda }\left[\nicefrac{3}{(\beta +1)^2}-\nicefrac{2}{(\beta +2)^2}\right]$, $\sym{cum}_{\beta\lambda}^{(\lambda)} =\nicefrac{n}{\lambda^2}\left(\nicefrac{3}{\beta+1}-\nicefrac{2}{\beta+2}\right)$, $\sym{cum}_{\lambda\lambda}^{(\beta)} =\nicefrac{2 n}{\lambda ^2}\left[\nicefrac{18}{(\beta +3)^2}-\nicefrac{8}{(\beta +4)^2}-\nicefrac{9}{(\beta +2)^2}\right]$ and $\sym{cum}_{\lambda\lambda}^{(\lambda)} =\nicefrac{2 n}{\lambda ^3}\left(1-\nicefrac{18}{\beta+2}+\nicefrac{36}{\beta+3}-\nicefrac{16}{\beta+4}\right)$.
  \end{coro}

  \Cref{prop:ECR:TD} presents the expected value of the third derivatives of the \ab{llf}, which can be obtained using essentially the same procedures employed in \cref{prop:ECR:FEI}.
  \begin{prop}\label{prop:ECR:TD}
    The expected value of the third derivatives of the \ab{llf} (third order cumulants of the \ab{llf}) in \eqref{eq:ECR:log-likelihood} with respect to its parameters are expressed by
    \begin{subequations}
      \begingroup
      \allowdisplaybreaks
      \begin{align}
        \sym{cum}_{\beta\beta\beta} & = \frac{2n}{\beta ^3},\label{eq:ECR:EVTDbbb} \\
        \sym{cum}_{\beta\beta\lambda} & = 0,\label{eq:ECR:EVTDbbl} \\
        \sym{cum}_{\beta\lambda\lambda} & = \frac{n}{\lambda^2}\left(\frac{9}{\beta +1}-\frac{28}{\beta +2}+\frac{27}{\beta +3}-\frac{8}{\beta +4}\right)\label{eq:ECR:EVTDbll} \intertext{and}
        \sym{cum}_{\lambda\lambda\lambda} & = \frac{2n}{\lambda ^3}\left(1-\frac{81}{\beta +2}+\frac{378}{\beta +3}-\frac{606}{\beta +4}+\frac{405}{\beta +5}-\frac{96}{\beta +6}\right). \label{eq:ECR:EVTDlll}
      \end{align}
      \endgroup
    \end{subequations}
  \end{prop}

  Let $ \hat{\theta}_i $, for $ i=1,\ldots,p $, be the \ab{MLE} for $\theta_i$. 
  \citet{CoxSnell-GeneralDefinitionResiduals-1968} proposed a method to improve \abp{MLE} in terms of the cumulants.
  This procedure was revisited by \citet{CordeiroKlein-Biascorrectionin-1994}, which established that the bias of  $ \hat{\theta}_i $ assumes the following formulation:
  \begin{equation}\label{eq:ECR:cox-snell-bias}
    \sym{bias}(\hat{\theta}_i)=\sym{E}(\hat{\theta}_i)-\theta_i=\sum_{r,s,t}\sym{cum}^{\theta_i,\theta_r}\sym{cum}^{\theta_s,\theta_t}\left(\sym{cum}_{\theta_r\theta_s}^{(\theta_t)}-\frac{1}{2}\sym{cum}_{\theta_r\theta_s\theta_t}\right) + \sum_{r=2}^\infty O(n^{-r}).
  \end{equation}
  Thus, the \ab{CS-MLE}, $ \widetilde{\theta}_i $, is defined by
  \begin{equation}\label{eq:ECR:cox-snell-correction}
    \widetilde{\theta}_i=\hat{\theta}_i-\widehat{\sym{bias}}(\hat{\theta}_i),
  \end{equation}
  where $ \widehat{\sym{bias}}(\hat{\theta}_i) $ represents the second order bias of $ \hat{\theta}_i $ evaluated in $ \hat{\theta}_i $. 
  \citet{CoxSnell-GeneralDefinitionResiduals-1968} showed that
  \begin{enumerate*}
    \item $ \sym{E}\left[\widehat{\sym{bias}}(\hat{\theta}_i)\right]=O(n^{-2}) $,
    \item $ \sym{E}(\hat{\theta}_i)=O(n^{-1}) $ and
    \item $ \sym{E}(\widetilde{\theta}_i)=\theta_i+O(n^{-2}) $.
  \end{enumerate*}
  Therefore, the bias of $ \widetilde{\theta}_i $ has order of $ n^{-2} $.
  Thus, the \abp{CS-MLE} are expected to possess better asymptotic behavior than \abp{MLE}.

  The \ab{MLE} second order biases are listed in \cref{prop:ECR:cox-snell-fisrtorderbias}.
  As expressed in \eqref{eq:ECR:cox-snell-correction}, they are the key to obtain the \abp{CS-MLE}.
  \begin{prop}\label{prop:ECR:cox-snell-fisrtorderbias}
    Let $X$ be an \ab{ECR} \ab{rv} with vector parameter $ \sym{PV}=(\beta,\lambda)^{\top} $. 
    The second order biases of the \ab{ECR} \abp{MLE} evaluated in $ \hat{\sym{PV}} $ are expressed as
    \begin{subequations}
      \begingroup
      \allowdisplaybreaks
      \begin{align*}
        \widehat{\sym{bias}}(\hat{\beta})=&\frac{1}{n}\left[\hat{\beta}^3
        +13 \hat{\beta} ^2+122 \hat{\beta} +380 -\frac{699840}{19321 (\hat{\beta} +5)}+\frac{96000}{361 (\hat{\beta} +6)}\right.\nonumber\\
        &+\frac{432 \left(4085783 \hat{\beta} ^2-8192586 \hat{\beta} -40352456\right)}{2641 \left(\hat{\beta} ^3-7 \hat{\beta} ^2+10 \hat{\beta} +72\right)^2}\nonumber\\
        &\left.-\frac{12 \left(70740551 \hat{\beta} ^2+3809213278 \hat{\beta} -35831044156\right)}{6974881 \left(\hat{\beta} ^3-7 \hat{\beta} ^2+10 \hat{\beta} +72\right)}\right]\intertext{and}
        \widehat{\sym{bias}}(\hat{\lambda})=&\frac{\hat{\lambda}}{n}\left[8 \hat{\beta}+86+\frac{49 }{270 \hat{\beta} }-\frac{1679616 }{96605 (\hat{\beta} +5)}+\frac{80000 }{1083 (\hat{\beta} +6)}\right.\nonumber\\
                                            &-\frac{8 \left(84037561 \hat{\beta} ^2 +21509105 \hat{\beta}  -393761162 \right)}{7923 \left(\hat{\beta} ^3-7 \hat{\beta} ^2+10 \hat{\beta} +72\right)^2}\nonumber\\
                                            &\left.+\frac{356431397749 \hat{\beta} ^2 -158970444943 \hat{\beta}  -4636191041858 }{376643574 \left(\hat{\beta} ^3-7 \hat{\beta} ^2+10 \hat{\beta} +72\right)}\right].
      \end{align*}
      \endgroup
    \end{subequations}
  \end{prop}

\begin{figure}[htb!]
    \centering
    \subfigure{
      \includegraphics[width=0.45\linewidth,trim=0 4in 0 0,clip]{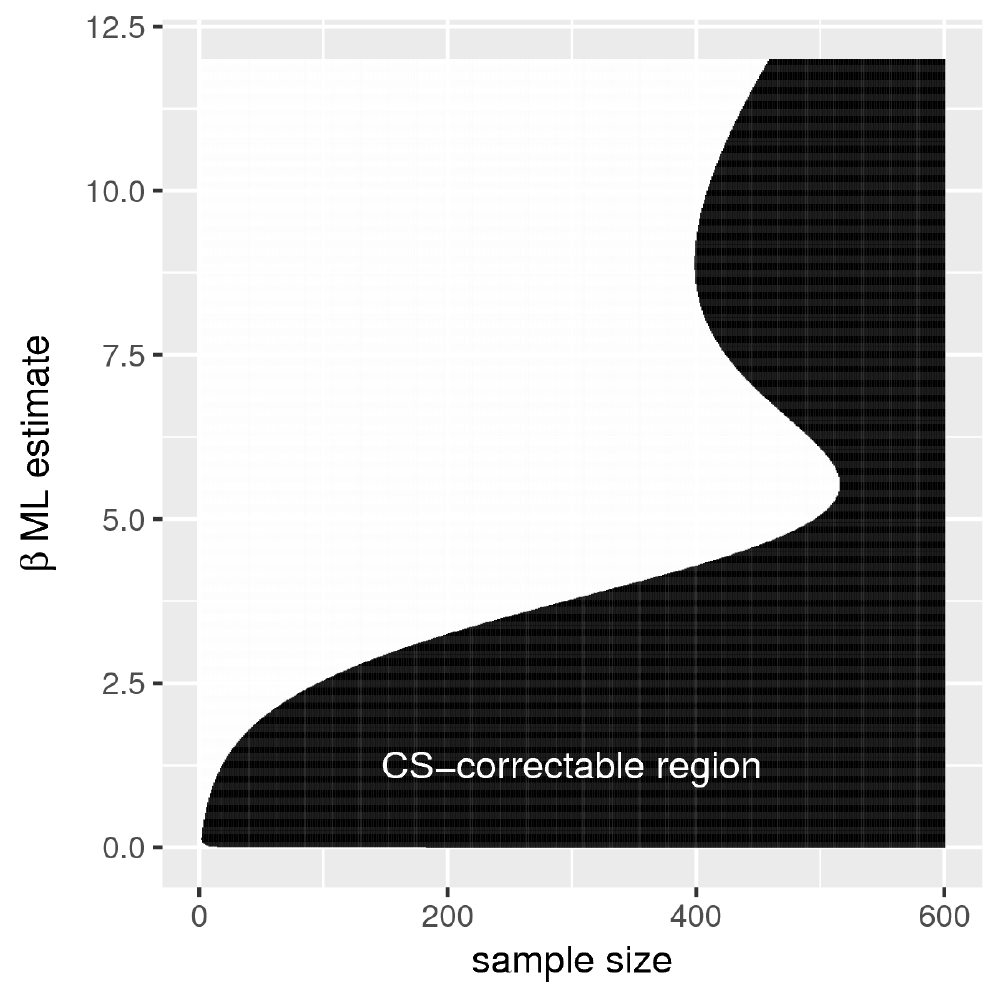}
    }
    \subfigure{
      \includegraphics[width=0.45\linewidth,trim=0 0 0 4in,clip]{Fig3.pdf}
    }
    \caption{CS-correctable region for \as{ECR} \as{ML} estimates.}
    \label{fig:ECR:cs-correctableregion}
  \end{figure}

  Unfortunately, for certain sample sizes and $ \hat{\beta} $ values, the \ab{ECR} \ab{CS-MLE} outcomes can be outside of the parametric space.
  Thus it is interesting to outline a map that identifies correctable regions, on which an constrained optimization would be acceptable.
  \Cref{fig:ECR:cs-correctableregion}  displays portions of the \ab{ECR} CS-correctable region in terms of sample size and $ \hat{\beta} $. 
  It is noticeable that this region do not depend of $ \hat{\lambda} $, as indicated by the $\widehat{\sym{bias}}(\widehat{\beta})$.

  The \ab{ECR} \abp{CS-MLE} are obtained using the \abp{MLE} second order biases expressed in \eqref{eq:ECR:cox-snell-correction}. 
  The \ab{ECR} \abp{CS-MLE} when $ \beta $ or $ \lambda $ are known can be obtained by \cref{prop:ECR:cox-snell-fisrtorderbiasknown}:
  \begin{prop}\label{prop:ECR:cox-snell-fisrtorderbiasknown}
    Let an \ab{ECR} \ab{rv} with parameters
    \begin{enumerate}
      \item $ \beta $ and $ \lambda=\lambda_0 $.
        The second order bias of the $ \beta $ \ab{MLE} evaluated in $ \hat{\beta} $ is expressed as
        \begin{equation*}
          \widehat{\sym{bias}}(\hat{\beta}|\lambda=\lambda_0)= \frac{\hat{\beta}}{n}.
        \end{equation*}
      \item $ \beta = \beta_0 $ and $ \lambda $.
        The second order bias of the $ \lambda $ \ab{MLE} evaluated in $ \hat{\lambda} $ is expressed as
        \begin{equation*}
          \widehat{\sym{bias}}(\hat{\lambda}|\beta=\beta_0)= \frac{\hat{\lambda}}{n}\frac{(\beta_0+2)(\beta_0+3)(\beta_0+4)}{\beta_0(\beta_0+5)(\beta_0+6)} \left[\frac{\beta_0^4+24 \beta_0^3+216 \beta_0^2+761 \beta_0 +294}{\left(\beta_0^2+11 \beta_0 +36\right)^2}\right].
        \end{equation*}
    \end{enumerate}
  \end{prop}

  The \ab{CR} \ab{CS-MLE} for $ \lambda $ is indicated in \cref{coro:CR:cox-snell-fisrtorderbias}.
  \begin{coro}\label{coro:CR:cox-snell-fisrtorderbias}
    Let a \ab{CR} \ab{rv} with parameter $ \lambda $.
    The second order bias of the $ \lambda $ \ab{MLE} evaluated in $ \hat{\lambda} $ is expressed as	$ \widehat{\sym{bias}}_{\as{CR}}(\hat{\lambda})= \nicefrac{45\hat{\lambda}}{56n} $.	
  \end{coro}

  \subsection{\Alp{PBE}}
  For data coming from a distribution having closed-form \ab{cdf}, a natural way to obtain unknown parameters is to determine the argument which minimizes the square distance between theoretical and sample percentiles.
  \citet{MurthyXieJiang-WeibullModels-WileySeriesInProbabilityAndStatistics-2004} discussed this method for the Weibull distribution, while \citet{GuptaKundu-Generalizedexponentialdistributiondifferentmethodofestimations-2001} studied the \ab{EE} case.

  Let $ X_{1:n}<X_{2:n}<\cdots<X_{n:n}$ be the order statistics obtained from a set of $ n $ \ab{iid} \abp{rv} following the \ab{ECR} distribution with vector parameter $ \sym{PV}=(\beta,\lambda)^{\top} $, the \ab{pbef} at $ \sym{PV} $ is given by
  \begin{equation}\label{eq:ECR:pbef}
    \sym{pbefS}(\sym{PV})=\sum_{i=1}^{n}\left[\frac{\lambda}{1-p_i^{\nicefrac{1}{\beta}}}\sqrt{\left(2-p_i^{\nicefrac{1}{\beta}}\right)p_i^{\nicefrac{1}{\beta}}}-x_{i:n}\right]^2.
  \end{equation}
  In this case, $p_i$ is the \ab{cdf} at the $i$th order statistic, $ \sym{cdfS}(x_{i:n}) $.
  \citet[p.~63]{MurthyXieJiang-WeibullModels-WileySeriesInProbabilityAndStatistics-2004} describe several ways to define sample $ p_i $.
  We use the \emph{mean rank} given by $ p_i=\nicefrac{i}{( n+1 )} $.

  The partial derivatives of the \ab{pbef} is provided by
  \begin{subequations}
    \begingroup
    \allowdisplaybreaks
    \begin{align}
      \frac{\partial}{\partial\beta}\sym{pbefS}(\sym{PV}) & = \frac{2\lambda}{\beta^2}
      [\lambda\hspace{-4ex}
      \underbrace{T_6(\beta,\bm{x})}_{\displaystyle\sum_{i=1}^{n}\frac{p_i^{\nicefrac{1}{\beta}} \log p_i}{\left(1-p_i^{\nicefrac{1}{\beta}}\right)^3}}
      \hspace{-4ex}-\hspace{-9ex}\overbrace{T_7(\beta,\bm{x})}^{\displaystyle\sum_{i=1}^{n}\frac{x_{i:n}\log p_i}{\left(1-p_i^{\nicefrac{1}{\beta}}\right)^2} \sqrt{\frac{p_i^{\nicefrac{1}{\beta}}}{2-p_i^{\nicefrac{1}{\beta}}}}}
      \hspace{-9ex}]
      \label{eq:ECR:dpbefbeta}
      \intertext{and}
      \frac{\partial}{\partial\lambda}\sym{pbefS}(\sym{PV}) & = 2
      [\hspace{-9ex}\underbrace{T_8(\beta,\bm{x})}_{-\displaystyle\sum_{i=1}^{n}\frac{x_{i:n}\sqrt{\left(2-p_i^{\nicefrac{1}{\beta}}\right)p_i^{\nicefrac{1}{\beta}}}}{1-p_i^{\nicefrac{1}{\beta}}}}
      \hspace{-9ex}-\lambda 
      \hspace{-6ex}
      \overbrace{T_9(\beta,\bm{x})}^{\displaystyle n-\sum_{i=1}^{n}\frac{1}{\left(1-p_i^{\nicefrac{1}{\beta}}\right)^2}}
      \hspace{-6ex}
      ]\label{eq:ECR:dpbeflambda}.
    \end{align}
    \endgroup
  \end{subequations}

  The \ab{PB} estimates, $\breve{\boldsymbol{\theta}}=(\breve{\beta},\breve{\lambda})^{\top}$, can be obtained numerically by solving the nonlinear equations system 
  $
  \nicefrac{\partial\sym{pbefS}(\sym{PV})}{\partial\beta}=\nicefrac{\partial\sym{pbefS}(\sym{PV})}{\partial\lambda}=0
  $.  
  In general, \abp{PBE} must be obtained numerically, similarly to \abp{MLE}.  
  Note that, from \cref{eq:ECR:dpbefbeta,eq:ECR:dpbeflambda}, it is possible to obtain two semi-closed \abp{PBE} for $\lambda$.
  For $\nicefrac{\partial\sym{pbefS}(\sym{PV})}{\partial\beta}=0$ and $\nicefrac{\partial\sym{pbefS}(\sym{PV})}{\partial\lambda}=0$ and a given $\beta$, the \abp{PBE} for $\lambda$ are expressed by, respectively:
  $\breve{\lambda}_1(\beta)=\nicefrac{T_7(\beta,\bm{x})}{T_6(\beta,\bm{x})}$
  and
  $\breve{\lambda}_2(\beta)=\nicefrac{T_8(\beta,\bm{x})}{T_9(\beta,\bm{x})}$.  
  Thus we can set the \ab{PB} estimate for $ \beta $ as solution the equation $ \breve{\lambda}_1(\beta)=\breve{\lambda}_2(\beta)  $ or, equivalently,
  $T_6(\beta,\bm{x})T_8(\beta,\bm{x})=T_7(\beta,\bm{x})T_9(\beta,\bm{x})$.
  By replacing $ \lambda $ by $ \breve{\lambda}_1(\beta) $ or $ \breve{\lambda}_2(\beta) $ in \eqref{eq:ECR:pbef}, we obtain the \abp{ppbef}:
  \begingroup
  \allowdisplaybreaks
  \begin{align*}
    \sym{pbefS}_{\lambda_1}(\beta)&= \sum_{i=1}^{n}\left[\frac{T_8(\beta,\bm{x})}{T_9(\beta,\bm{x})\left(1-p_i^{\nicefrac{1}{\beta}}\right)}\sqrt{\left(2-p_i^{\nicefrac{1}{\beta}}\right)p_i^{\nicefrac{1}{\beta}}}-x_{i:n}\right]^2\intertext{and}
    \sym{pbefS}_{\lambda_2}(\beta)&= \sum_{i=1}^{n}\left[\frac{T_7(\beta,\bm{x})}{T_6(\beta,\bm{x})\left(1-p_i^{\nicefrac{1}{\beta}}\right)}\sqrt{\left(2-p_i^{\nicefrac{1}{\beta}}\right)p_i^{\nicefrac{1}{\beta}}}-x_{i:n}\right]^2.
  \end{align*}
  \endgroup
  We can also obtain the \ab{PB} estimates for $ \beta $ by minimizing one of the \abp{ppbef} with respect to $ \beta $.

  \section{Simulation studies}\label{sec:ECR:SimulationStudy}

  In this section, we assess the performance of the three proposed estimators by means of two simulation studies.
  First we investigate the convergence of the bias-corrected \abp{MLE}.
  Subsequently, the performance of \ab{MLE}, \ab{CS-MLE} and \ab{PBE} is quantified adopting relative biases and \abp{SSD} as comparison criteria.
  To that end, we used an Intel\textsuperscript \textregistered\ Core\texttrademark\ i7-4500U CPU at 1:80 GHz processor with base x64, Ubuntu 16.04.2 operating system, and the \texttt{R}-3.3.3 \citep{RFoundationforStatisticalComputing-RLanguageand-2008} computational platform.

  \subsection{Convergence study}\label{sec:ECR:ConvergenceStudy}

  In this study, we adopt sample sizes $ n\in\{5, 10, \ldots, 250\} $ and $\sym{PV}=(0.5,0.6)$.
  For each pair $ (\boldsymbol{\theta},n) $, we generate 1,000 Monte Carlo replications using the \ab{itsm}. 
  For each Monte Carlo replica, we obtained both bias and \ab{SSD} and, subsequently, their Monte Carlo averages are calculated.
  \begin{figure}[htb!]
    \centering
    \subfigure[bias for $ \beta=0.5 $]{
      \includegraphics[width=0.33\textwidth,trim=0 298 310 45,clip]{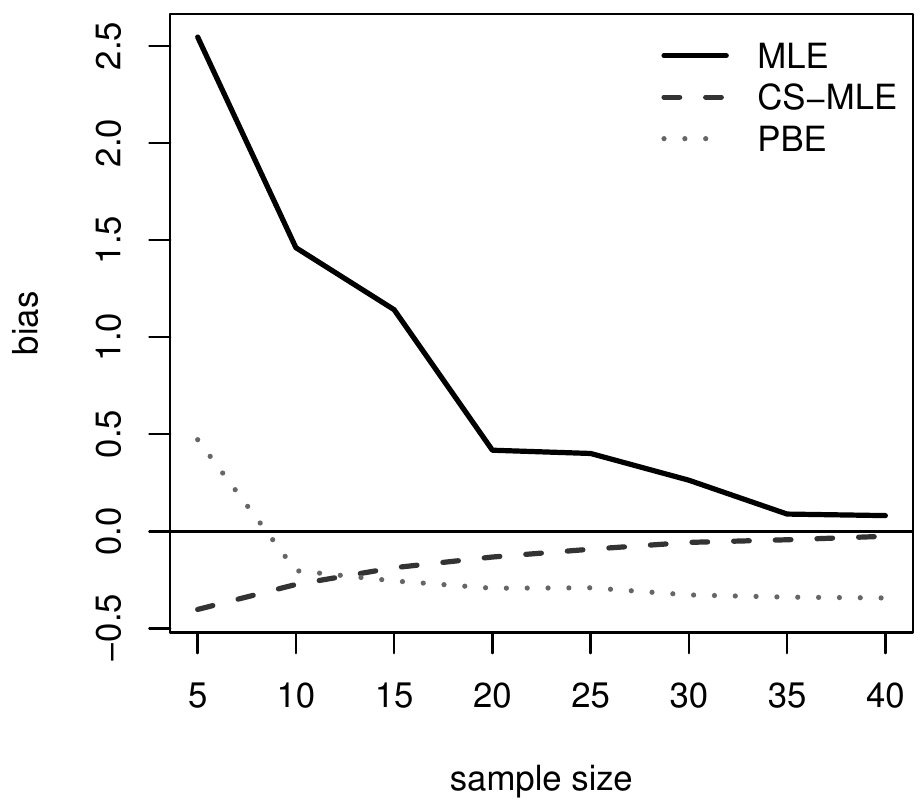}
    }
    \subfigure[bias for $ \lambda=0.6 $]{
      \includegraphics[width=0.33\textwidth,trim=290 298 20 45,clip]{Fig4.pdf}
    }
    \subfigure[\as{SSD} for $ \beta=0.5 $]{
      \includegraphics[width=0.33\textwidth,trim=0 10 310 333,clip]{Fig4.pdf}
    }
    \subfigure[\as{SSD} for $ \lambda=0.6 $]{
      \includegraphics[width=0.33\textwidth,trim=290 10 20 333,clip]{Fig4.pdf}
    }
    \caption{$ \operatorname{\as{ECR}}( 0.5,0.6) $ mean bias and \as{SSD} for various sample sizes.}
    \label{fig:ECR:ConvStudyB0.5L0.6}
  \end{figure}

  \Cref{fig:ECR:ConvStudyB0.5L0.6} exhibits bias and \ab{SSD} values for several sample sizes.
  It is noticeable the $\hat{\beta}$ bias outperforms that due to $\breve{\beta}$ for smaller size samples.
  Under CS-correctable region, $\tilde{\beta}$ presents the smallest bias for the majority of cases.
  With respect to \ab{SSD}, $\breve{\beta}$ tends to assume the smallest value at small samples; while, $\tilde{\beta}$ has the best asymptotic features. 

  In general, the $\hat{\lambda}$ bias is the smallest, while the $\breve{\lambda}$ one is the highest. 
  The $\hat{\lambda}$ and $\tilde{\lambda}$ \abp{SSD} present similar values, while $\breve{\lambda}$ is the highest one.


  \subsection{Relative biases and \asp{SSD} study}\label{sec:ECR:relativeBiasStudy}

\begin{figure}[htb!]
    \centering
    \subfigure[$ \beta $ estimate relative biases]{
      \includegraphics[width=0.33\textwidth,trim=0 298 310 45,clip]{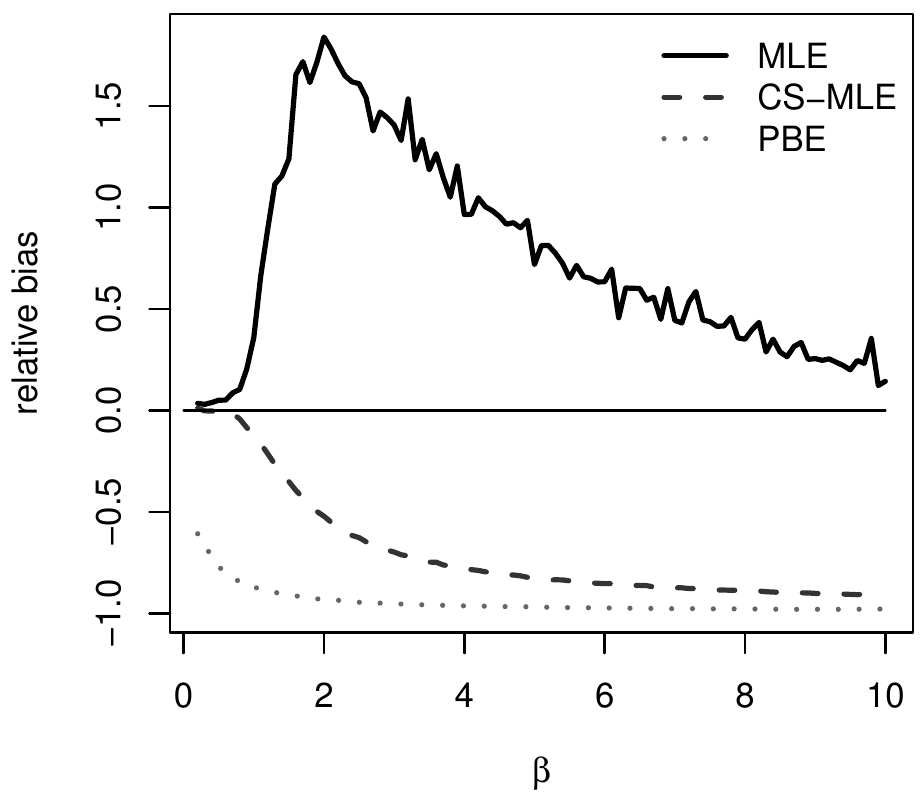}
      \label{fig:ECR:ComparisonBiasSSDn100-betabias}
    }
    \subfigure[$ \lambda $ estimate relative biases]{
      \includegraphics[width=0.33\textwidth,trim=290 298 20 45,clip]{Fig5.pdf}
      \label{fig:ECR:ComparisonBiasSSDn100-lambdabias}
    }
    \subfigure[$ \beta $ estimate relative \abp{SSD}]{
      \includegraphics[width=0.33\textwidth,trim=0 10 310 333,clip]{Fig5.pdf}
      \label{fig:ECR:ComparisonBiasSSDn100-betassd}
    }
    \subfigure[$ \lambda $ estimate relative \abp{SSD}]{
      \includegraphics[width=0.33\textwidth,trim=290 10 20 333,clip]{Fig5.pdf}
      \label{fig:ECR:ComparisonBiasSSDn100-lambdassd}
    }
    \caption{Comparison between relatives bias and \as{SSD} of \as{ECR} estimates for $ \lambda=1 $ and sample size 100.}
    \label{fig:ECR:ComparisonBiasSSDn100}
  \end{figure}

\begin{figure}[htb!]
    \centering
    \subfigure[\asp{MLE} of $ \beta $]{
      \includegraphics[width=0.31\textwidth,trim=0 6in 4in 0,clip]{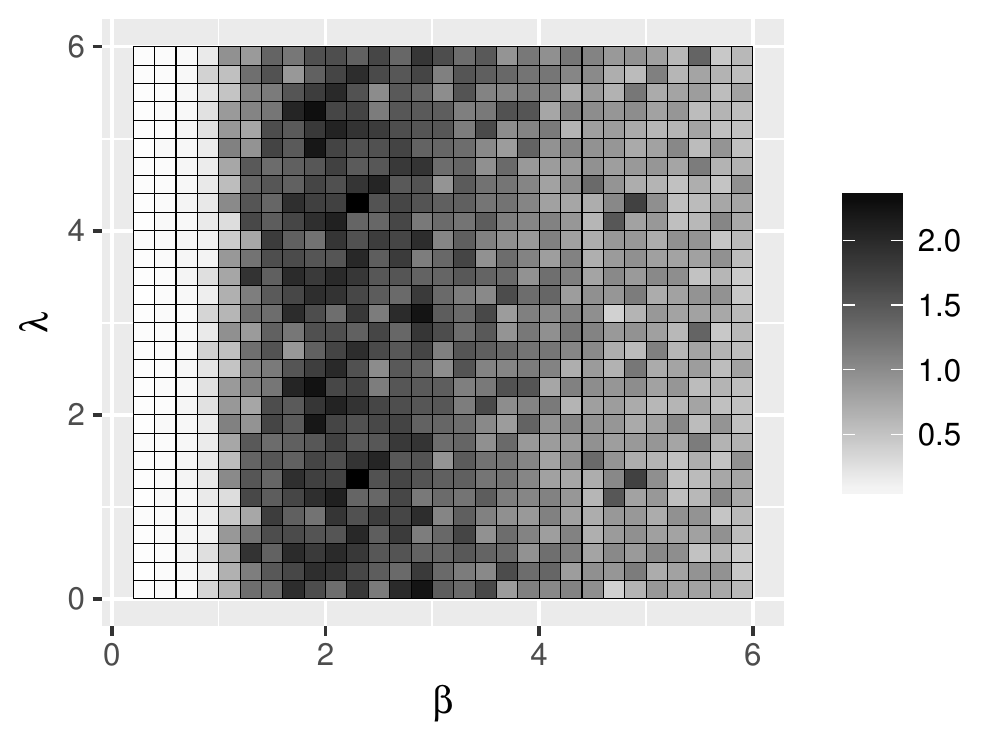}
    }
    \subfigure[\asp{CS-MLE} of $ \beta $]{
      \includegraphics[width=0.31\textwidth,trim=0 3in 4in 3in,clip]{Fig6.pdf}
    }
    \subfigure[\asp{PBE} of $ \beta $]{
      \includegraphics[width=0.31\textwidth,trim=0 0 4in 6in,clip]{Fig6.pdf}
    }
    \subfigure[\asp{MLE} of $ \lambda $]{
      \includegraphics[width=0.31\textwidth,trim=4in 6in 0 0,clip]{Fig6.pdf}
    }
    \subfigure[\asp{CS-MLE} of $ \lambda $]{
      \includegraphics[width=0.31\textwidth,trim=4in 3in 0 3in,clip]{Fig6.pdf}
    }
    \subfigure[\asp{PBE} of $ \lambda $]{
      \includegraphics[width=0.31\textwidth,trim=0 3in 0 0,clip]{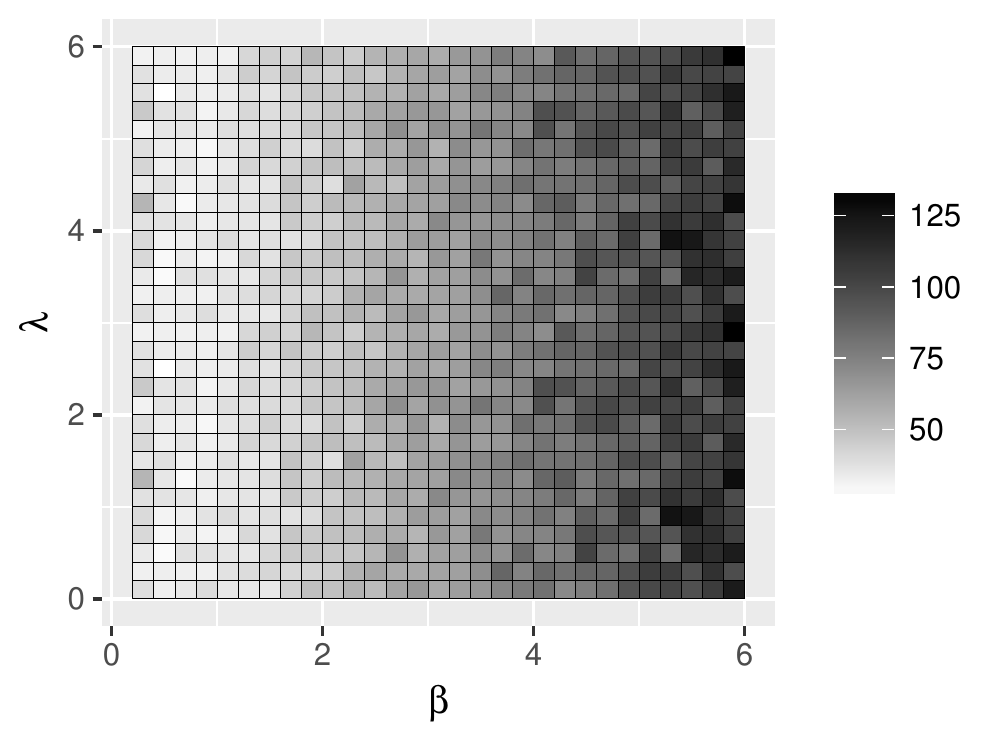}
    }
    \caption{\as{ECR} relative biases estimates at the sample size 100.}
    \label{fig:ECR:SimulationStudyBias}
  \end{figure}

  In this study, we regard to two sample sizes, $n_1=10$ and $n_2=100$, and the parameter vectors $\boldsymbol{\theta}=(\beta,\lambda)$ such that $\beta,\lambda\in\{0.1, 0.2, \ldots,10\}$.
  For each pair $(\boldsymbol{\theta},n_i)$, $i = 1, 2$, we perform 1,000 Monte Carlo replications.
  For each random sample, we calculated both bias and \ab{SSD}.
  Finally, we compute the average of biases and \abp{SSD}.
  Results are depicted in \crefrange{fig:ECR:SimulationStudyBias}{fig:ECR:ComparisonBiasSSDn100}.

\begin{figure}[htb!]
    \centering
    \subfigure[\asp{MLE} of $ \beta $]{
      \includegraphics[width=0.31\textwidth,trim=0 6in 4in 0,clip]{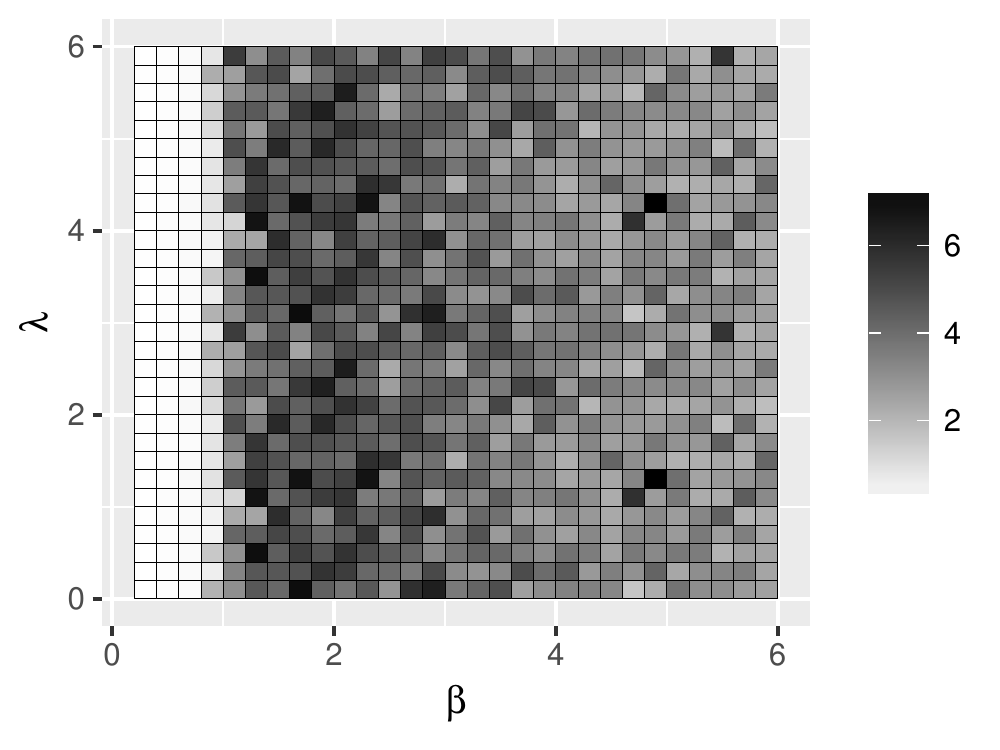}
    }
    \subfigure[\asp{CS-MLE} of $ \beta $]{
      \includegraphics[width=0.31\textwidth,trim=0 3in 4in 3in,clip]{Fig7.pdf}
    }
    \subfigure[\asp{PBE} of $ \beta $]{
      \includegraphics[width=0.31\textwidth,trim=0 0 4in 6in,clip]{Fig7.pdf}
    }
    \subfigure[\asp{MLE} of $ \lambda $]{
      \includegraphics[width=0.31\textwidth,trim=4in 6in 0 0,clip]{Fig7.pdf}
    }
    \subfigure[\asp{CS-MLE} of $ \lambda $]{
      \includegraphics[width=0.31\textwidth,trim=4in 3in 0 3in,clip]{Fig7.pdf}
    }
    \subfigure[\asp{PBE} of $ \lambda $]{
      \includegraphics[width=0.31\textwidth,trim=4in 0 0 6in,clip]{Fig7.pdf}
    }
    \caption{\as{ECR} relative \as{SSD} at the estimates for sample size 100.}
    \label{fig:ECR:SimulationStudyStdErr}
  \end{figure}

  In \cref{fig:ECR:ComparisonBiasSSDn100}, we can observe that the $\hat{\beta}$ relative bias and \ab{SSD} decrease for when $\beta$ increases, while the relative bias and \ab{SSD} of $\tilde{\beta}$ and $\breve{\beta}$ do not have this property. 
  Small values of $\beta$, the smallest relative biases and \abp{SSD} for $\beta$ estimates are obtained by $\tilde{\beta}$ and $\breve{\beta}$.
  In general, the smallest biases and \abp{SSD} for $\lambda$ estimates are obtained by $\hat{\lambda}$ and $\tilde{\lambda}$.
  In particular, for smallest values of $\beta$, the $\tilde{\lambda}$ has the smallest biases and \abp{SSD}.
  In \cref{fig:ECR:SimulationStudyBias,fig:ECR:SimulationStudyStdErr} biases and \abp{SSD} do not appear to suffer influence of $\lambda$ value, \as{ie}, biases and \abp{SSD} are to be only function of $\beta$ values.
  Figures also reveal that \abp{CS-MLE} are desirable for only small values of $\beta$.

  \section{Application}\label{sec:ECR:Application}

  An application to real data is done to illustrate the \ab{ECR} potentiality.
  To that end, we use an uncensored real dataset previously discussed by \citet{CrowleyHu-CovarianceAnalysisHeart-1977}, which consists of 66 patient survival times with respect to a Stanford heart transplant list.
  We consider only patients that died in the followup time and were not submitted to prior bypass surgery.
  The data are presented in \cref{table:ECR:CrowleyHudataset}.
  \begin{table}[htb!]
    \centering
    \caption{\citeauthor{CrowleyHu-CovarianceAnalysisHeart-1977}'s dataset}
    \begin{tabular}{rrrrrrrrrrr}
      \hline
      1 &   1 &   2 &   2 &   2 &   4 &   4 &   5 &   5 &    7 &    8 \\
      11 &  15 &  15 &  15 &  16 &  17 &  20 &  20 &  27 &   29 &   31 \\
      34 &  35 &  36 &  38 &  39 &  42 &  44 &  49 &  50 &   52 &   57 \\
      60 &  65 &  67 &  67 &  68 &  71 &  71 &  76 &  77 &   79 &   80 \\
      84 &  89 &  95 &  99 & 101 & 109 & 148 & 152 & 187 &  206 &  218 \\
      262 & 284 & 284 & 307 & 333 & 339 & 674 & 732 & 851 & 1031 & 1386 \\
      \hline
    \end{tabular}
    \label{table:ECR:CrowleyHudataset}
  \end{table}
  \Cref{table:ECR:DescriptiveStatistics} gives a descriptive summary, which anticipates data are right-skewed.
  \begin{table}[htb!]
    \centering
    \caption[Descriptive Statistics for \citeauthor{CrowleyHu-CovarianceAnalysisHeart-1977}'s dataset]{Descriptive Statistics}
    \begin{tabular}{lrrrrrrr}
      \hline
      Dataset                                             &   Mean & Median & Variance &    Skew. &    Kurt. & Min. & Max. \\
      \hline
      \citeauthor{CrowleyHu-CovarianceAnalysisHeart-1977} & 143.70 &  58.50 & 64506.42 & 3.104648 & 13.00242 &    1 & 1386 \\
      \hline
    \end{tabular}
    \label{table:ECR:DescriptiveStatistics}
  \end{table}

\begin{figure}[htb!]
    \centering
    \subfigure[\as{TTT}-plot]{
      \includegraphics[width=0.30\textwidth,trim=0 10 0 55,clip]{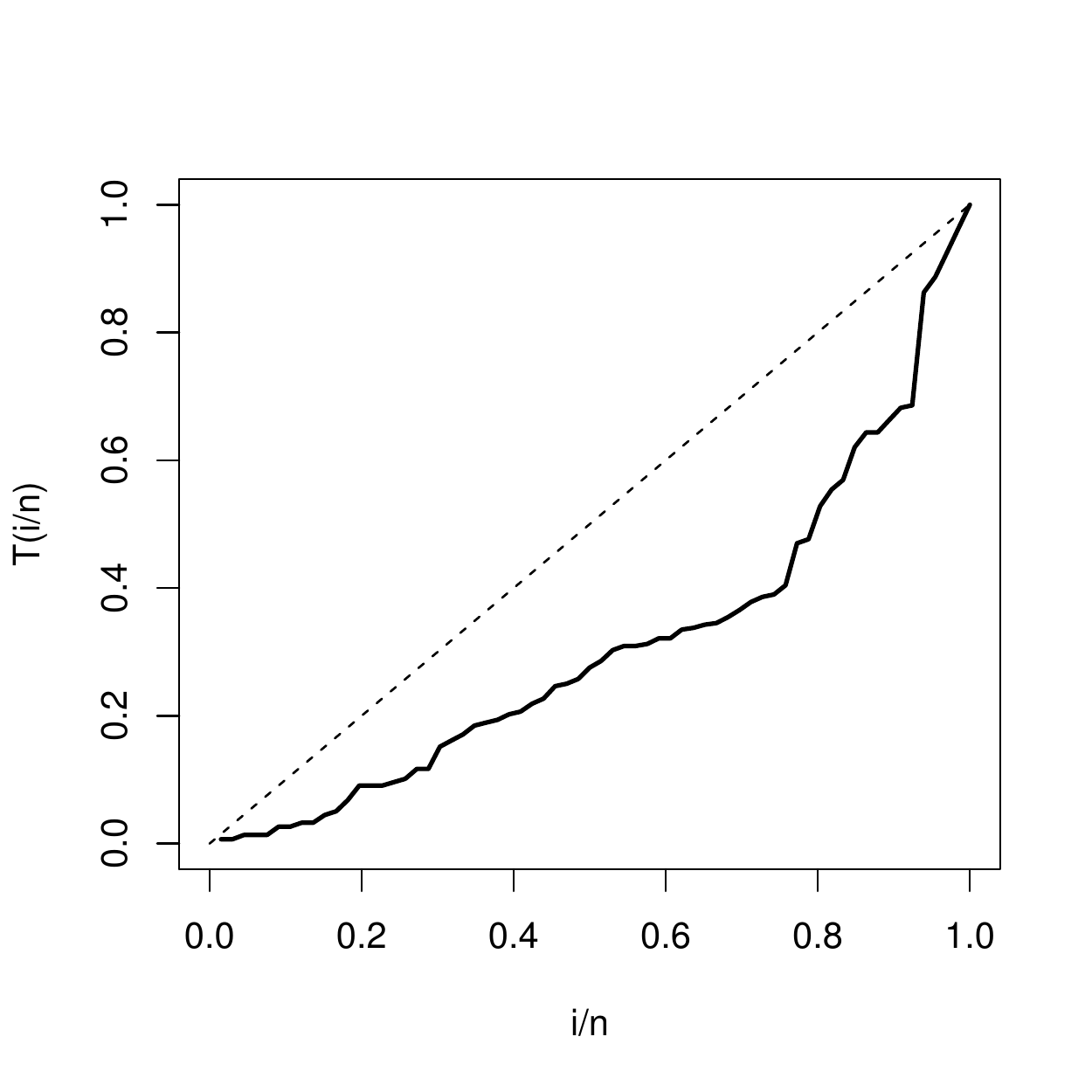}
      \label{fig:ECR:TTT}
    }
    \subfigure[log-likelihood surface]{
      \includegraphics[width=0.30\textwidth,trim=30 60 30 70,clip]{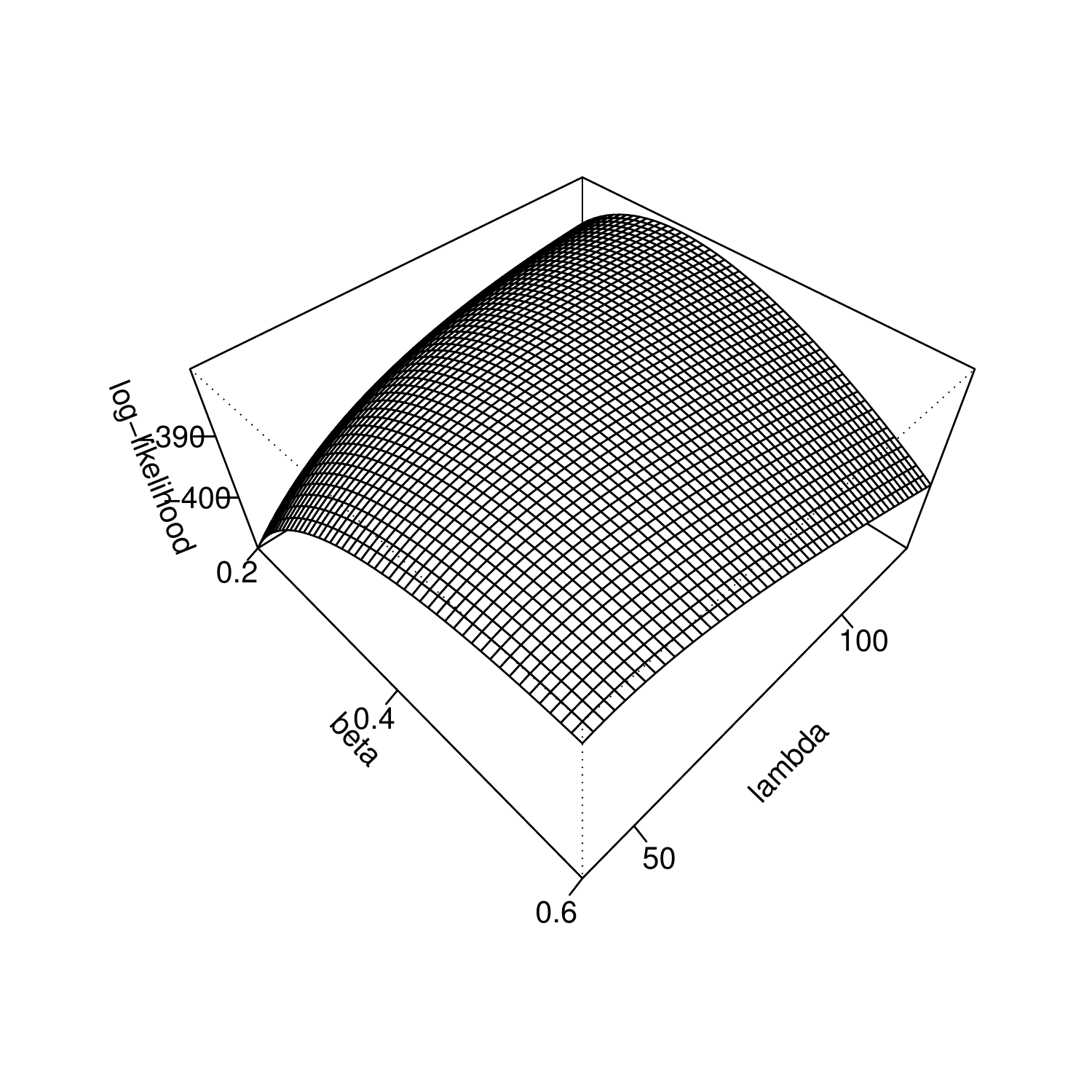}
      \label{fig:ECR:log-likelihoodsurface}
    }
    \subfigure[log-likelihood contour]{
      \includegraphics[width=0.30\textwidth,trim=0 20 10 50,clip]{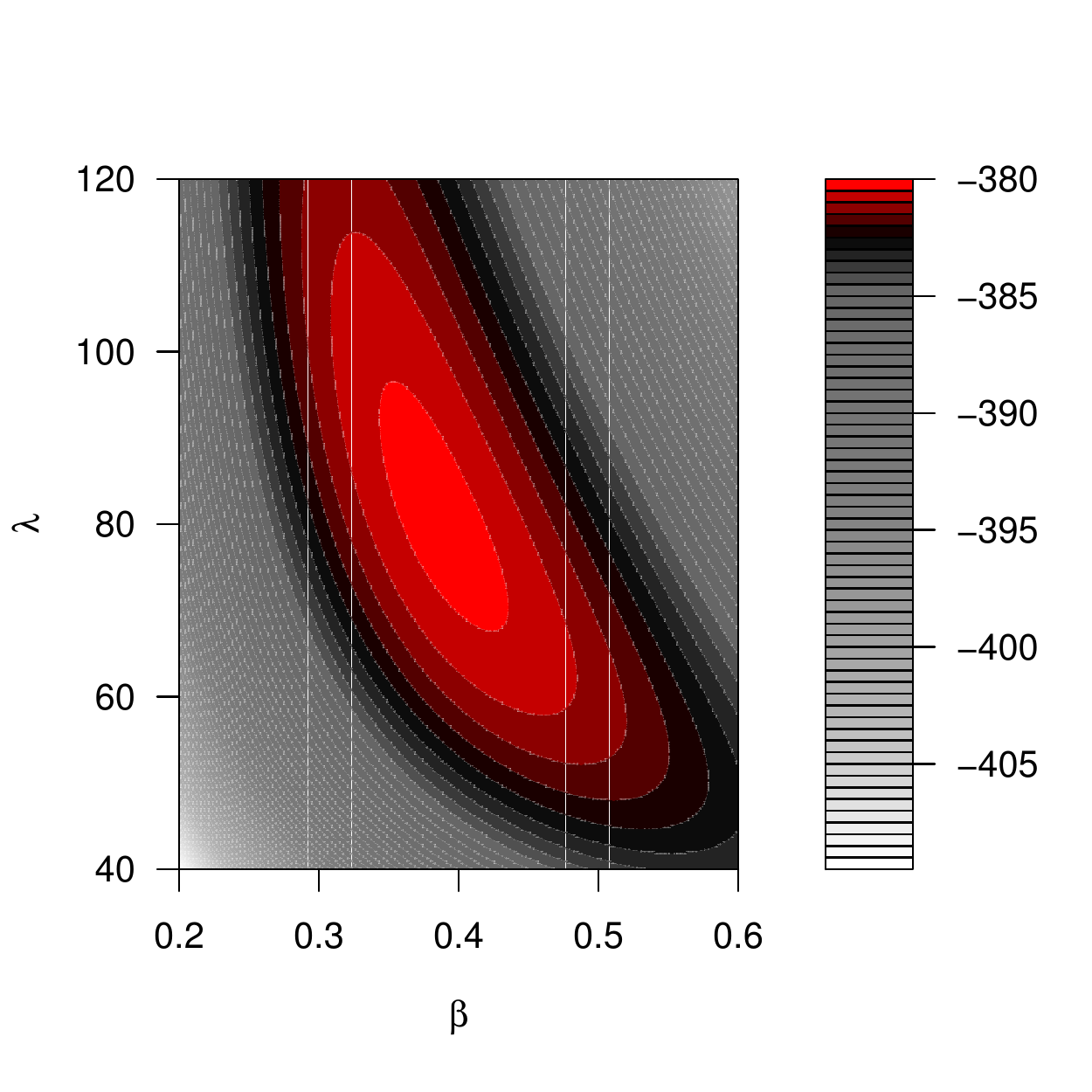}
      \label{fig:ECR:log-likelihoodcontour}
    }
    \caption{\as{TTT} and log-likelihood plots of \citeauthor{CrowleyHu-CovarianceAnalysisHeart-1977}'s dataset.}
    \label{fig:ECR:TTTandqqplots}
  \end{figure}

  In many applications, empirical \ab{hrf} shape can indicate a particular model.
  The plot of \ab{TTT} \citep{Aarset-HowtoIdentify-1987} can be useful in this sense.
  The \ab{TTT} plot is obtained by measuring $G(\nicefrac{r}{n})= \nicefrac{( \sum_{i=1}^{r}y_{i:n}+(n-r)y_{r:n})}{\sum_{i=1}^{n}y_{i:n}}$ against $\nicefrac{r}{n}$ for $r=1,\ldots,n$.
  The under study \ab{TTT} plot is presented in \cref{fig:ECR:TTT}.
  Results indicate a decreasing \ab{hrf}, pattern which is covered by the \ab{ECR} \ab{hrf}.
  Beyond, \cref{fig:ECR:log-likelihoodsurface,fig:ECR:log-likelihoodcontour} present the log-likelihood surface and its contour curves, respectively, in order to anticipate the \ab{ECR} optimization framework at under study data.
  The concavity of the surface $\sym{llfS}(\beta,\lambda)$ indicates an extreme local maximum over an elliptical region centred around $(\beta_0,\lambda_0)=(0.4,80)$.

  To confirm the \ab{ECR} relevance among biparametric models, we compare it with other twenty distributions:  
  \begin{itemize}
    \item The \ab{GG} model with shape parameters $\mu,q\in\sym{Re}$ and scale parameter $\sigma>0$ as given in \citet{Prentice-LogGammaModel-1974};
    \item The \ab{SM} model with shape parameters $ c>0 $ and $k>0$ and scale parameter $\lambda>0$ obtained from \cref{eq:compoundWeibull} assuming $a=0$;
    \item The \ab{EHC} model with shape parameter $ a>0 $ and scale parameter $ \phi>0 $ as given in \citet{CordeiroLemonte-BetaHalfCauchy-2011};
    \item The \ab{LL} model with shape parameter $a\ge1$ and scale parameter $b>0$ as given \citet[p.~223]{KleiberKotz-StatisticalSizeDistributionsinEconomicsandActuarialSciences-2003};
    \item The \ab{LN} model with shape parameters $ \mu>0 $ and $\sigma>0$ as given in \citet[p.~107]{KleiberKotz-StatisticalSizeDistributionsinEconomicsandActuarialSciences-2003};
    \item The Weibull model with shape parameter $ a>0 $ and scale parameter $ \beta>0$ as given in \citet[p.~174]{KleiberKotz-StatisticalSizeDistributionsinEconomicsandActuarialSciences-2003};
    \item The \ab{LE} model with shape parameter $ \theta>0 $ and scale parameter $\lambda>0$ as given in \citet{BhatiMalikVaman-LindleyExponentialdistribution-2015};
    \item The \ab{GHN} model with shape parameter $ \alpha>0 $ and scale parameter $ \theta>0 $ as given in \citet{CoorayAnanda-GeneralizationHalfNormal-2008};
    \item The Fréchet (inverse Weibull) model with shape parameter $ \lambda>0 $ and scale parameter $ \sigma>0$ as given in \citet[p.~294]{Bury-StatisticalDistributionsinEngineering-1999};
    \item The Chen model with shape parameter $ \beta>0 $ and scale parameter $ \lambda>0 $ as given in \citet{Chen-newtwoparameter-2000};
    \item The two-parameter \ab{BSd} model with shape parameters $\alpha>0$ and $\beta>0$ as given in \citet{BirnbaumSaunders-newfamilylife-1969};
    \item The gamma model with shape parameter $ p>0 $ and scale parameter $ b>0 $ as given in \citet[p.~148]{KleiberKotz-StatisticalSizeDistributionsinEconomicsandActuarialSciences-2003};
    \item The \ab{EE} model with shape parameter $ \alpha>0 $ and scale parameter $ \lambda>0 $ as given in \citet{GuptaKundu-ExponentiatedExponentialFamily-2001};
    \item The \ab{ELi} model with shape parameter $ \alpha>0 $ and scale parameter $ \lambda>0 $ as given in \citet{NadarajahBakouchTahmasbi-generalizedLindleydistribution-2011};
    \item The \ab{IG} model with shape parameter $ \alpha>0 $ and scale parameter $ \beta>0 $ as given in \citet{LeemisMcQueston-UnivariateDistributionRelationships-2008};
    \item The \ab{BX} model with shape parameter $ \alpha>0 $ and scale parameter $\lambda>0$ as given in \citet{KunduRaqab-GeneralizedRayleighdistribution-2005};
    \item The Gompertz model with shape parameter $ a>0 $ and scale parameter $ b>0 $ as given in \citet{Lenart-Gompertzdistributionand-MfdsF-2012};
    \item The Wald (inverse normal) model with shape parameter $ \mu>0 $ and scale parameter $\lambda>0$ as given in \citet{MichaelSchucanyHaas-GeneratingRandomVariates-1976};
    \item The \ab{FW} model with shape parameter $\alpha>0$ and scale parameter $\beta>0$ as given in \citet{BebbingtonLaiZitikis-flexibleWeibullextension-2007};
    \item The \ab{BXII} model with shape parameters $ c>0 $ and $k>0$ obtained from \cref{eq:compoundWeibull} assuming $a=0$ and $\lambda=1$;
    \item The \ab{CR} model with scale parameter $\lambda>0$.
  \end{itemize}

\begin{figure}[htb!]
    \centering
    \subfigure[\asp{pdf} vs. histogram]{
      \includegraphics[width=0.45\textwidth,trim=0 10 0 55,clip]{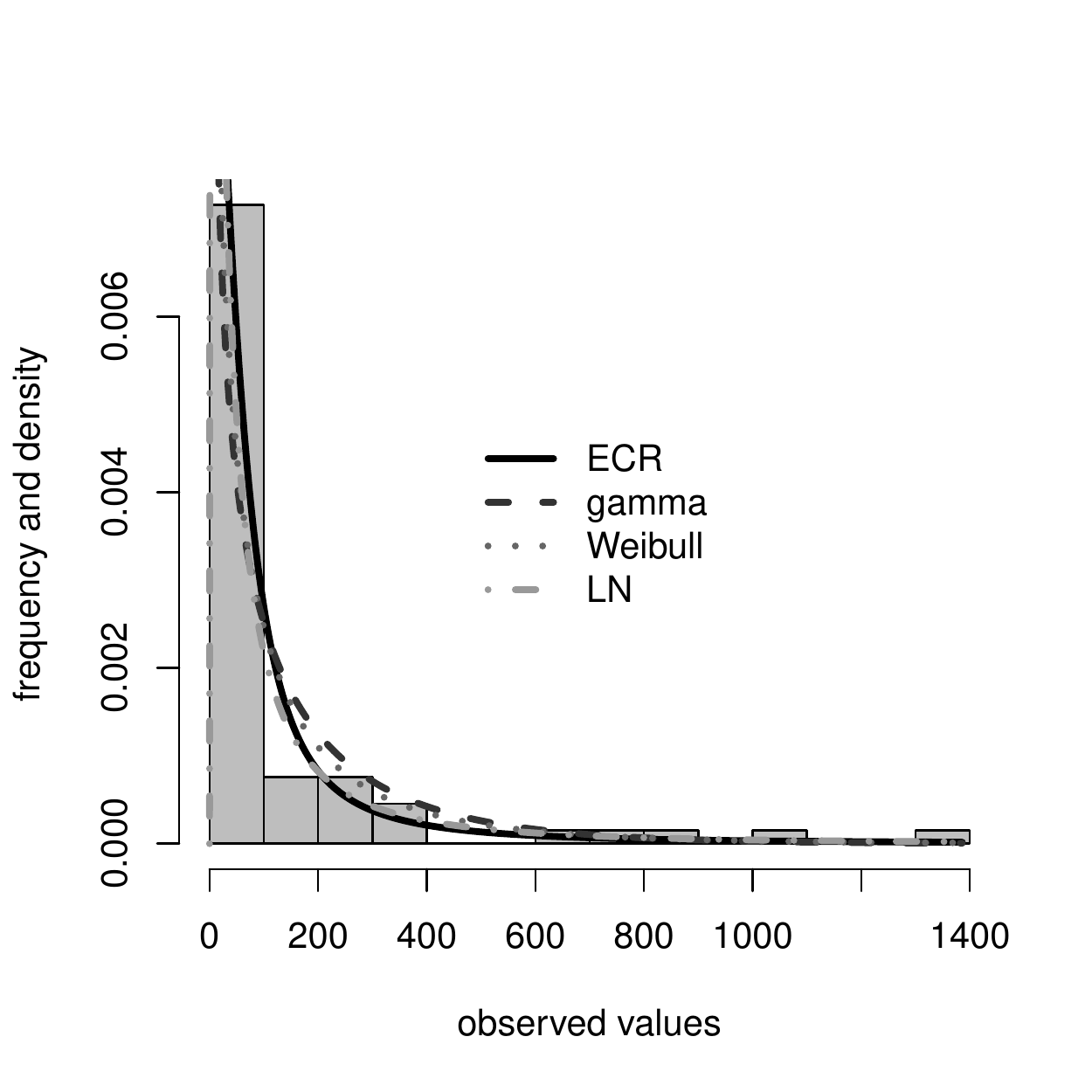}
      \label{fig:ECR:ApplicationPDF}
    }
    \subfigure[\asp{cdf} vs. \asp{ecdf}]{
      \includegraphics[width=0.45\textwidth,trim=0 10 0 55,clip]{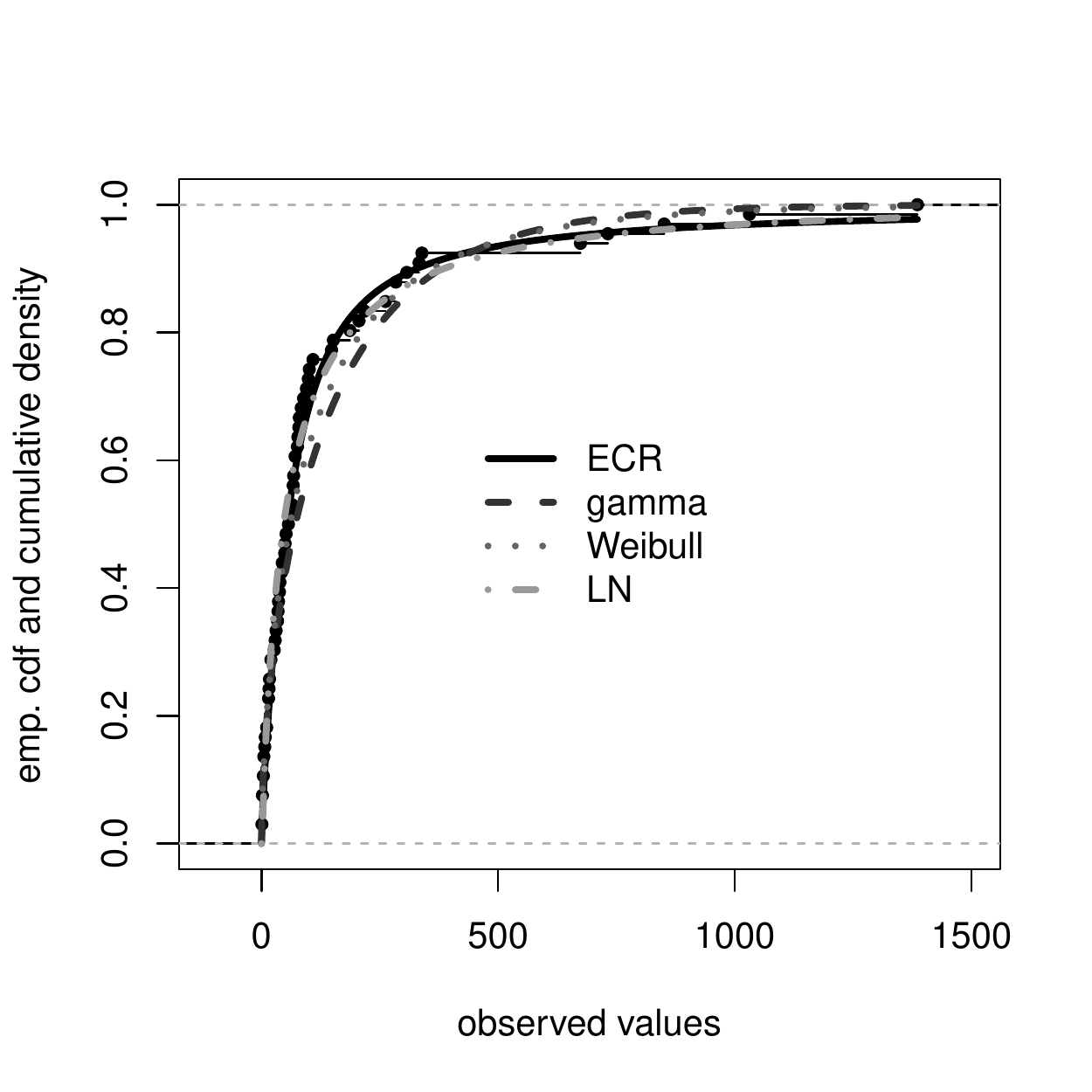}
      \label{fig:ECR:ApplicationCDF}
    }
    \caption{Plots of \as{ECR}, gamma, Weibull and \as{LN} models.}
  \end{figure}

  \cref{fig:ECR:ApplicationPDF} displays fitted and empirical densities.
  For better exhibition, we regard only the \ab{ECR}, \ab{LN}, Weibull and gamma models which seem to provide closer fits to the histogram.
  The empirical and fitted \abp{cdf} of these models are displayed in \cref{fig:ECR:ApplicationCDF}.
  From this plot, we note that the \ab{ECR} model provides a good fit.
  The \ab{ECR} qq-plot is depicted in \cref{fig:ECR:qqplot} and it reveals that \ab{ECR} fitted quantiles are very close to sample quantiles.
  \Cref{table:ECR:EstimativesApplication} presents the \ab{ML} estimates and their std. errors for fitted models, indicating there is not non-significant fits from their respective asymptotic confidence intervals.
  \begin{table}[h!]
    \centering
    \caption{Estimates for \citeauthor{CrowleyHu-CovarianceAnalysisHeart-1977}'s dataset}
    \begin{small}
      \begin{tabular}{lrlrlrl}
        \hline
        model                      &  estimate & (std. error) &           &            &           &            \\
        \hline
        \ab{GG}$(\mu,\sigma,q)$    & 4.14959   & (0.32326) & 1.58916   & (0.14812)  & 0.36973   & (0.15358)  \\
        \ab{SM}$(c,k,\lambda)$     & 1.89529   & (0.05652) & 0.91101   & (0.26596)  & 129.90197 & (10.56896) \\
        \ab{ECR}$(\beta,\lambda)$  & 0.38669   & (0.04002) & 80.68399  & (11.70962) &           & \\
        \ab{EHC}$(a,\phi)$         & 0.78024   & (0.09915) & 70.31135  & (4.33446)  &           & \\
        \ab{LL}$(a,b)$             & 1.08065   & (0.11220) & 50.33793  & (4.19446)  &           & \\
        \ab{LN}$(\mu,\sigma)$      & 3.84913   & (0.20212) & 1.64286   & (0.14354)  &           & \\
        Weibull$(a,\beta)$         & 0.66923   & (0.05851) & 104.10812 & (4.20483)  &           & \\
        \ab{LE}$(\theta,\lambda)$  & 0.79155   & (0.09721) & 0.00372   & (0.00152)  &           & \\
        \ab{GHN}$(\alpha,\theta)$  & 0.49637   & (0.04215) & 142.78850 & (4.20515)  &           & \\
        Fréchet$(\lambda,\sigma)$  & 0.58458   & (0.05110) & 20.11157  & (4.45158)  &           & \\
        Chen$(\beta,\lambda)$      & 0.21043   & (0.01224) & 0.06867   & (0.01512)  &           & \\
        \ab{BSd}$(\alpha,\beta)$   & 2.26636   & (0.19821) & 37.62267  & (4.21313)  &           & \\
        gamma$(p,b)$               & 0.55830   & (0.04101) & 257.40580 & (40.42949) &           & \\
        \ab{EE}$(\alpha,\lambda)$  & 0.55028   & (0.08222) & 0.00449   & (0.00112)  &           & \\
        \ab{ELi}$(\alpha,\lambda)$ & 0.27401   & (0.04215) & 0.00577   & (0.00107)  &           & \\
        \ab{IG}$(\alpha,\beta)$    & 0.45398   & (0.06615) & 5.14991   & (1.23915)  &           & \\
        \ab{BX}$(\beta,\lambda)$   & 0.19605   & (0.02615) & 0.00175   & (0.00025)  &           & \\
        Gompertz$(a,b)$            & 0.00688   & (0.00041) & 0.00011   & (0.00009)  &           & \\
        Wald$(\mu,\lambda)$        & 143.70599 & (4.84315) & 12.31591  & (2.42212)  &           & \\
        \ab{FW}$(\alpha,\beta)$    & 0.00147   & (0.00057) & 8.34328   & (1.41122)  &           & \\
        \ab{BXII}$(c,k)$           & 0.01528   & (0.00125) & 16.99507  & (2.56261)  &           & \\
        \ab{CR}$(\lambda)$         & 24.49100  & (0.44792) &           &            &           & \\
        \hline
      \end{tabular}
    \end{small}
    \label{table:ECR:EstimativesApplication}
  \end{table}

  To compare quantitatively discussed models, \cref{table:ECR:GoodnessOfFitApplication} elects the \ab{GoF} statistics of the fitted models.
  We use the following \ab{GoF} statistics:
  \begin{itemize}
    \item Criteria under \abp{cdf}:
      \begin{itemize}
        \item \ab{W*};
        \item \ab{A*};
        \item \ab{KS}.
      \end{itemize}
    \item Criteria under \abp{pdf}:
      \begin{itemize}
        \item \ab{AIC};
        \item \ab{BIC};
        \item \ab{CAIC};
        \item \ab{HQIC}.
      \end{itemize}
  \end{itemize}
  Smaller \ab{GoF} values under \abp{cdf} are associated with better fits, while the first \abp{GoF} class may indicate superiority relations in nested models.
  \begin{table}[h!]
    \centering
    \caption{\Al{GoF} tests for \citeauthor{CrowleyHu-CovarianceAnalysisHeart-1977}'s dataset}
    \begin{tabular}{lrrrrrrr}
      \hline
      model                      &        \ab{W*} &        \ab{A*} &        \ab{KS} &         \ab{AIC} &        \ab{CAIC} &         \ab{BIC} &        \ab{HQIC} \\
      \hline
      \ab{GG}$(\mu,\sigma,q)$    &          0.069 &          0.399 &          0.079 &          765.532 &          765.919 &          772.101 &          768.127 \\
      \ab{SM}$(c,k,\lambda)$     &          0.052 &          0.327 &          0.072 &          766.012 &          766.399 &          772.581 &          768.608 \\
      \ab{ECR}$(\beta,\lambda)$  & \textbf{0.039} & \textbf{0.286} & \textbf{0.057} &          \textbf{764.612} &          \textbf{764.803} &          \textbf{768.992} &          \textbf{766.343} \\
      \ab{EHC}$(a,\phi)$         &          0.041 &          0.295 &          0.060 &          764.627 &          764.817 &          769.006 &          766.357 \\
      \ab{LL}$(a,b)$             &          0.074 &          0.454 &          0.063 &          765.481 &          765.672 &          769.861 &          767.212 \\
      \ab{LN}$(\mu,\sigma)$      &          0.102 &          0.579 &          0.089 &          764.915 &          765.105 &          769.294 &          766.645 \\
      Weibull$(a,\beta)$         &          0.114 &          0.689 &          0.118 &          767.444 &          767.635 &          771.824 &          769.175 \\
      \ab{LE}$(\theta,\lambda)$  &          0.140 &          0.838 &          0.139 &          768.735 &          768.926 &          773.115 &          770.466 \\
      \ab{GHN}$(\alpha,\theta)$  &          0.191 &          1.150 &          0.142 &          773.212 &          773.403 &          777.592 &          774.943 \\
      Fréchet$(\lambda,\sigma)$  &          0.356 &          2.036 &          0.146 &          780.374 &          780.564 &          784.753 &          782.104 \\
      Chen$(\beta,\lambda)$      &          0.295 &          1.766 &          0.151 &          782.206 &          782.397 &          786.585 &          783.937 \\
      \ab{BSd}$(\alpha,\beta)$   &          0.229 &          1.229 &          0.155 &          770.982 &          771.173 &          775.362 &          772.713 \\
      gamma$(p,b)$               &          0.195 &          1.172 &          0.159 &          772.658 &          772.849 &          777.037 &          774.389 \\
      \ab{EE}$(\alpha,\lambda)$  &          0.210 &          1.261 &          0.168 &          773.709 &          773.900 &          778.088 &          775.440 \\
      \ab{ELi}$(\alpha,\lambda)$ &          0.258 &          1.562 &          0.185 &          778.603 &          778.793 &          782.982 &          780.333 \\
      \ab{IG}$(\alpha,\beta)$    &          0.563 &          3.183 &          0.210 &          792.664 &          792.854 &          797.043 &          794.394 \\
      \ab{BX}$(\beta,\lambda)$   &          0.386 &          2.274 &          0.238 &          788.552 &          788.742 &          792.931 &          790.282 \\
      Gompertz$(a,b)$            &          0.215 &          1.291 &          0.242 &          793.824 &          794.014 &          798.203 &          795.554 \\
      Wald$(\mu,\lambda)$        &          0.425 &          2.383 &          0.261 &          787.709 &          787.899 &          792.088 &          789.439 \\
      \ab{FW}$(\alpha,\beta)$    &          0.931 &          4.691 &          0.262 &          816.638 &          816.829 &          821.018 &          818.369 \\
      \ab{BXII}$(c,k)$           &          0.676 &          3.828 &          0.323 &          824.816 &          825.007 &          829.196 &          826.547 \\
      \ab{CR}$(\lambda)$         &          0.213 &          1.254 &          0.132 &          785.023 &          785.085 &          787.212 &          785.888 \\
      \hline
    \end{tabular}
    \label{table:ECR:GoodnessOfFitApplication}
  \end{table}
  From figures of merit, the proposed \ab{ECR} distribution presents smallest \ab{GoF} values.
  To test $\mathcal{H}_ 0:\beta=1$ (or $\sym{cdfS}_{\operatorname{\ab{CR}}}=\sym{cdfS}_{\operatorname{\ab{ECR}}}$), we employ the likelihood ratio statistic, which yields a $p$-value $<10^{-5}$, indicating there is statistical difference between \ab{ECR} and \ab{CR} models. 
  In summary, the \ab{ECR} model may be a good alternative for describing heavy-tailed positive real data.

  \section{Concluding remarks}\label{sec:ECR:Concludingremarks}

  We deepen a discussion about an extended Cauchy-Rayleigh distribution, particular case of some known models like \al{HTR} \citep{NikiasShao-SignalProcessingwithAlpha-1995}, \al{GFP} \citep{Zandonatti-DistribuzionideParetoGeneralizzate-UniversityofTrento-2001}, \al{EBXII} \citep{Al-HussainiAhsanullah-ExponentiatedDistributions-AtlantisStudiesinProbabilityandStatistics5-2015}, \al{KBXII} \citep{ParanaibaOrtegaCordeiroPascoa-KumaraswamyBurrXII-2013} and \al{McBXII} \citep{GomesSilvaCordeiro-TwoExtendedBurr-2015}.
  We refer to it as \af{ECR} distribution.
  Some of its mathematical properties are derived and discussed:
  Closed-form expressions for mode and probability weighted, log-, incomplete and order statistic moments.
  The \as{ECR} model obeys the property of regularly varying at infinity and can take decreasing, decreasing-increasing-decreasing and upside-down bathtub-shape \alp{hrf}, confirming its usefulness to describe lifetime data.
  We provide procedures to estimate the \as{ECR} parameters through original and bias-corrected \alp{MLE} and a \al{PB} method.
  A simulation study to assess proposed estimators is performed.
  Each procedure reveals advantages over specific parameter points and sample sizes.
  The \as{ECR} usefulness is illustrated through an application to real data.
  Results indicate that our proposal can furnish better performance than classical lifetime models like gamma, \al{BSd}, Weibull and \al{LN}.
  We hope that the broached model may attract wider applications for modeling positive real data.

  \section{Appendix}

  \begin{proof}[Proof of \cref{prop:ECR:limitspdf}]
    Note that
    \begingroup
    \allowdisplaybreaks
    \begin{align*}
      \lim\limits_{x\to 0^+}f(x) & = \beta\lambda \lim\limits_{x\to 0^+}\frac{x}{(\lambda^2+x^2)^{\nicefrac{3}{2}}}\left(1- \frac{\lambda}{\sqrt{\lambda^2+x^2}}\right)^{\beta-1}\frac{\left(1+ \frac{\lambda}{\sqrt{\lambda^2+x^2}}\right)^{\beta-1}}{\left(1+ \frac{\lambda}{\sqrt{\lambda^2+x^2}}\right)^{\beta-1}}\\
                                 & = \beta\lambda \lim\limits_{x\to 0^+}\frac{x^{2\beta-1}}{(\lambda^2+x^2)^{\nicefrac{1}{2}+\beta}}\left(1+ \frac{\lambda}{\sqrt{\lambda^2+x^2}}\right)^{1-\beta}.
    \end{align*}
    \endgroup
    Now the result in \eqref{eq:ECR:limitspdf} is obtained.	
  \end{proof}

  \begin{proof}[Proof of \cref{prop:ECR:rvi}]
    Let $ X $ be a \ab{ECR} \ab{rv} and $ c>0 $.
    \begingroup
    \allowdisplaybreaks
    \begin{align*}
      \lim_{x\to\infty}\frac{\sym{sfS}(cx)}{\sym{sfS}(x)} 
      & = \lim_{x\to\infty}\frac{1-\left(1-\frac{\lambda}{\sqrt{\lambda^2+(cx)^2}}\right)^\beta}{1-\left(1-\frac{\lambda}{\sqrt{\lambda^2+x^2}}\right)^\beta} \\
      & = c^2 \lim_{x\to\infty}\left\{\left[\frac{\lambda^2+x^2}{\lambda^2+(cx)^2}\right]^{\nicefrac{3}{2}}\,
    \frac{\left(1-\frac{\lambda}{\sqrt{\lambda^2+(cx)^2}}\right)^{\beta-1}}{\left(1-\frac{\lambda}{\sqrt{\lambda^2+x^2}}\right)^{\beta-1}}\right\}\quad\text{(L'Hôpital's rule)}\\
    & = c^2\left(\frac{1}{c^2}\right)^{\nicefrac{3}{2}}=\frac{1}{c}.
  \end{align*}
  \endgroup
\end{proof}

\begin{proof}[Proof of \cref{prop:ECR:ECRmode}]
  Note that
  \begingroup
  \allowdisplaybreaks
  \begin{equation}
    \frac{\mathrm{d} \log\left[f(x)\right]}{\mathrm{d}x} = \frac{\lambda ^2 \left[(\beta -1) \lambda +\beta  \sqrt{\lambda ^2+x^2}\,\right]-x^2 \left[(1-\beta)\lambda +2 \sqrt{\lambda ^2+x^2}\right]}{x \left(\lambda ^2+x^2\right)^{\nicefrac{3}{2}}}\label{eq:ECR:FSlogpdf}.
  \end{equation}
  \endgroup	
  By \cref{eq:ECR:FSlogpdf} the mode is obtained solving
  \begin{equation*}
    \lambda ^2 \left[(\beta -1) \lambda +\beta  \sqrt{\lambda ^2+x^2}\,\right]-x^2 \left[(1-\beta)\lambda +2 \sqrt{\lambda ^2+x^2}\right]=0.
  \end{equation*}	
  After some manipulations we obtain
  \begin{equation}\label{eq:ECR:mode:chareq1}
    \frac{2x^2-\lambda^2\beta}{\sqrt{\lambda^2+x^2}}=(\beta-1)\lambda.
  \end{equation}	
  For $x_0$ to be a solution of \cref{eq:ECR:mode:chareq1} it needs to satisfy the following conditions:
  \begin{subequations}
    \begin{align}
      \beta>1 &\iff x_0>\lambda\sqrt{\nicefrac{\beta}{2}}\label{eq:ECR:cec1}\\
      0<\beta<1 &\iff 0<x_0<\lambda\sqrt{\nicefrac{\beta}{2}}\label{eq:ECR:cec2},
    \end{align}
  \end{subequations}	
  the case $\beta=1$ is trivial.	
  Squaring both sides of \cref{eq:ECR:mode:chareq1} and simplifying we obtain a bi-quadratic equation
  \begin{equation}\label{eq:ECR:mode:QEx}
    4x^4-\left[(\beta+1)\lambda\right]^2 x^2+\lambda^4(2\beta-1)=0.
  \end{equation}	
  Assuming $y=x^2$ the last expression became
  \begin{equation}\label{eq:ECR:mode:QEy}
    4y^2-\left[(\beta+1)\lambda\right]^2 y+\lambda^4(2\beta-1)=0.
  \end{equation}	
  The \al{Delta} of this quadratic equation is $ \sym{Delta} = (\beta -1)^2 \left(\beta ^2+6 \beta +17\right) \lambda ^4 $, and its roots are given by
  $
  y_1=\nicefrac{\lambda^2}{8}\left[(\beta+1)^2+\sqrt{(\beta-1)^2(\beta^2+6\beta+17)}\right]
  $
  and 
  $
  y_2=\nicefrac{\lambda^2}{8}\left[(\beta+1)^2-\sqrt{(\beta-1)^2(\beta^2+6\beta+17)}\right]
  $.
  This solutions can be rewritten as
  \begin{align*}
    y_1^\star&=\frac{\lambda^2}{8}\left[(\beta+1)^2+(\beta-1)\sqrt{\beta^2+6\beta+17}\right]\\
    y_2^\star&=\frac{\lambda^2}{8}\left[(\beta+1)^2-(\beta-1)\sqrt{\beta^2+6\beta+17}\right].
  \end{align*}	
  Note that for $ \beta>1 $ we have $ y_1=y_1^\star $ and $ y_2=y_2^\star $, on the other hand if $ \beta<1 $ we obtain $ y_1=y_2^\star $ and $ y_2=y_1^\star $ and finally, when $\beta=1$ all solutions are equivalent.	
  It is straightforward to prove that
  $
  y_1^\star\ge 0\iff\beta\ge\nicefrac{1}{2}
  $ and 
  $
  y_2^\star\ge 0\iff\beta\ge 0
  $.
  Then \cref{eq:ECR:mode:QEx} has only two possible (real and non-negative) roots (for $\beta\ge\nicefrac{1}{2}$)
  \begin{subequations}
    \begin{align}
      x_1&=\frac{\lambda}{2\sqrt{2}}\sqrt{(\beta+1)^2+(\beta-1)\sqrt{\beta^2+6\beta+17}},\quad\beta\ge\nicefrac{1}{2},\label{eq:ECR:modesol}\\
      x_2&=\frac{\lambda}{2\sqrt{2}}\sqrt{(\beta+1)^2-(\beta-1)\sqrt{\beta^2+6\beta+17}},\quad\beta\ge 0.\label{eq:ECR:modesol2}
    \end{align}
  \end{subequations}
  \Cref{eq:ECR:modesol2} do not satisfy any of the conditions in \eqref{eq:ECR:cec1} and \eqref{eq:ECR:cec2}, then it is not a solution of \cref{eq:ECR:mode:chareq1}.
  On the contrary the solution in \eqref{eq:ECR:modesol} satisfy both \cref{eq:ECR:cec1} and 
  $
  \nicefrac{1}{2}\le\beta<1 \iff 0\le x_1<\lambda\sqrt{\nicefrac{\beta}{2}}
  $,
  which is (using \cref{coro:ECR:ECRmodellimits}) a restricted version of \cref{eq:ECR:cec2}.
  Then the unique solution of \cref{eq:ECR:mode:chareq1} is \cref{eq:ECR:modesol}.
  Note that the \ab{ECR} \ab{pdf} is not defined in zero and hence the \ab{ECR} mode is given by \cref{eq:ECR:ECRmode}.
\end{proof}

\begin{proof}[Proof of \cref{prop:ECR:auxprop1}]Note that
  $$\sym{Pr}(X<x)=\sym{Pr}(\lambda Z<x)=\sym{Pr}\left(Z< \frac{x}{\lambda}\right)=\left(1-\frac{1}{\sqrt{1-\left( \frac{x}{\lambda}\right)^2}}\right)^\beta=\left(1- \displaystyle\frac{\lambda}{\sqrt{\lambda^2+x^2}}\right)^\beta $$
\end{proof}

\begin{proof}[Proof of \cref{prop:ECR:pwms}]
  Let $X\sim\operatorname{\as{ECR}}(\beta,\lambda)$.
  Note that
  \begin{equation*}
    \sym{pwmS} = \beta\lambda\int_{0}^{\infty}x^r\frac{x}{\left(\lambda^2+x^2\right)^{\nicefrac{3}{2}}}\left(1-\frac{\lambda}{\sqrt{\lambda^2+x^2}}\right)^{(s+1)\beta-1}\left[1-\left(1-\frac{\lambda}{\sqrt{\lambda^2+x^2}}\right)^\beta\right]^k\,\mathrm{d}x.
  \end{equation*}
  We apply a simple change of variable assuming
  \begin{equation}\label{eq:ECR:cv1}
    u=1-\frac{\lambda}{\sqrt{\lambda^2+x^2}}\iff x=\lambda\frac{\sqrt{u(2-u)}}{1-u},
  \end{equation}	
  and then after some algebras we obtain
  \begin{equation*}
    \sym{pwmS} = \beta(\lambda\sqrt{2})^r\int_0^1\frac{u^{\nicefrac{r}{2}+(s+1)\beta-1}(1-u)^{-r}}{(1-\nicefrac{u}{2})^{\nicefrac{-r}{2}}}(1-u^\beta)^t\,\mathrm{d}u.
  \end{equation*}
  Using simple binomial expansion we can write $ (1-u^\beta)^t = \sum_{i=0}^{t}(-1)^i\binom{t}{i}i^{i\beta}  $.
  Then we can express
  \begin{equation}\label{eq:ECR:pwmstep1}
    \sym{pwmS} = \beta(\lambda\sqrt{2})^r\sum_{i=0}^{t}(-1)^i\binom{t}{i}\int_0^1\frac{u^{\nicefrac{r}{2}+(s+i+1)\beta-1}(1-u)^{-r}}{(1-\nicefrac{u}{2})^{\nicefrac{-r}{2}}}\,\mathrm{d}u.
  \end{equation}
  Consider the following result
  \begin{equation}\label{eq:ECR:2F1integral}
    \sym{2F1}(a,b;c;x)=\frac{1}{\sym{bf}(b,c-b)}\int_{0}^{1}\frac{t^{b-1}(1-t)^{c-b-1}}{(1-tx)^a}\,\mathrm{d}t,
  \end{equation}	
  showed in \citet[p.~4]{Bailey-Generalizedhypergeometricseries-1973}.
  The result is obtained applying \cref{eq:ECR:2F1integral} in \eqref{eq:ECR:pwmstep1}.
\end{proof}

\begin{proof}[Proof of \cref{coro:ECR:momentsCR}]
  Let $Z\sim\operatorname{\as{CR}}(1)$.
  Consider the following identity obtained by \citet[p.~11]{Bailey-Generalizedhypergeometricseries-1973}
  \begin{equation}\label{eq:ECR:hi1}
    \sym{2F1}(a,1-a;c;\nicefrac{1}{2})=\frac{\sym{gf}\left(\frac{c}{2}\right)\sym{gf}\left(\frac{1+c}{2}\right)}{\sym{gf}\left(\frac{a+c}{2}\right)\sym{gf}\left(\frac{1+c-a}{2}\right)}.
  \end{equation}	
  Assuming $\beta=1$ in \eqref{eq:ECR:momentsECR}, applying \cref{eq:ECR:hi1} and after some algebras we can set
  \begin{equation*}
    \sym{E}(Z^r) = \frac{2^{\nicefrac{r}{2}+1}}{\sqrt{\pi }}
    \frac{\sym{gf}\left(\frac{1}{2} \left(2-\frac{r}{2}\right)\right) \sym{gf} \left(\frac{1}{2} \left(3-\frac{r}{2}\right)\right)}{ \sym{gf} \left(\frac{2-r}{2}\right)}\,\sym{bf}\left(1-r,\frac{r}{2}+1\right) = \frac{\sym{gf}\left( \frac{1-r}{2}\right)\,\sym{gf}\left(1+ \frac{r}{2}\right)}{\sqrt{\pi}}.
  \end{equation*}	
  Then the results in \eqref{eq:momentsCR} follows.
\end{proof}

\begin{proof}[Proof of \cref{prop:ECR:ElogX}]
  Note that $ \sym{E}(\log Z) = \beta\int_{0}^{\infty}\!\log z \frac{z}{(1+z^2)^{\nicefrac{3}{2}}} \left\{1- \frac{1}{\sqrt{1+z^2}}\right\}^{\beta-1}\,\mathrm{d}z $. We apply a simple change of variable assuming $ u=1-\nicefrac{1}{\sqrt{1+z^2}}$ and after some algebras we can set
  \begin{equation*}
    \sym{E}(\log Z) = 
    \frac{\beta}{2}\int_{0}^{1}\!u^{\beta-1}\log u\,\mathrm{d}u + 
    \frac{\beta}{2}\int_{0}^{1}\!u^{\beta-1}\log (2-u)\,\mathrm{d}u - 
    \beta\int_{0}^{1}\!u^{\beta-1}\log (1-u)\,\mathrm{d}u.
  \end{equation*}	
  The first integral is determined transforming $ v=-\log u $. After this transform we obtain
  \begin{equation}\label{eq:ECR:intu}
    \int_{0}^{1}\!u^{\beta-1}\log u\,\mathrm{d}u=-\int_{0}^{\infty}\,v\exp(-\beta v)\,\mathrm{d}v=-\frac{1}{\beta^2}.
  \end{equation}	
  In order to determine the second integral we need to look upon the following known expansion $ \log(2-u)=\log 2-\sum_{n=1}^{\infty}\nicefrac{1}{n}\left(\nicefrac{u}{2}\right)^n  $.
  Then we can express the second integral as
  \begingroup
  \allowdisplaybreaks
  \begin{align}
    \int_{0}^{1}\!u^{\beta-1}\log (2-u)\,\mathrm{d}u & = (\log 2)\int_{0}^{1}\!u^{\beta-1}\,\mathrm{d}u-\sum_{n=1}^{\infty}\frac{1}{n2^n}\int_{0}^{1}u^{\beta+n-1}\,\mathrm{d}u \nonumber           \\
                                                     & = \frac{\log 2}{\beta} - \sum_{n=1}^{\infty}\frac{1}{n(n+\beta)2^n}                                                                = \frac{\log 2}{\beta} - \frac{1}{\beta}\sum_{n=1}^{\infty}\left[\frac{1}{n2^n}-\frac{1}{(n+\beta)2^n}\right]                      \nonumber \\
                                                     & = \frac{1}{\beta}\left[\sum_{n=0}^{\infty}\frac{1}{(n+\beta)2^n}-\frac{1}{\beta}\right] = \frac{1}{\beta}\left[\sym{lpf}\left(\frac{1}{2};1,\beta\right)-\frac{1}{\beta}\right].\label{eq:ECR:int2u}
  \end{align}
  \endgroup
  Consider another well known expansion $ \log(1-u)=-\sum_{n=1}^{\infty}\nicefrac{u^n}{n} $, then the last integral can be expressed as
  \begin{align}
    \int_{0}^{1}\!u^{\beta-1}\log(1-u) \,\mathrm{d}u & =-\sum_{n=1}^{\infty}\frac{1}{n}\int_{0}^{1}\!u^{\beta+n-1}\,\mathrm{d}u =-\sum_{n=1}^{\infty}\frac{1}{n(\beta+n)}                                \nonumber\\
                                                     & =-\frac{1}{\beta}\left[\sym{dgf}(1+\beta)+\sym{EM}\right]. \label{eq:ECR:int1u}
  \end{align}
  In the last step we use the identity obtained from \citet[eq.~6.3.16]{AbramowitzStegun-Handbookofmathematicalfunctions-1972}: $ \sym{dgf}(1+\beta) =  - \sym{EM} + \sum_{n=1}^{\infty}\nicefrac{\beta}{[ n(\beta+n) ]} $.
  Then we can set $ \sym{E}(\log Z)=\nicefrac{1}{2}\sym{lpf}\left(\nicefrac{1}{2};1,\beta\right)+\sym{dgf}(1+\beta)+\sym{EM}-\nicefrac{1}{\beta} $, and the result follows noting that $\sym{E} (\log X) =\log\lambda+\sym{E} (\log Z)$.
\end{proof}

\begin{proof}[Proof of \cref{prop:ECR:ims}]
  The \abp{im} for a \ab{rv} $ X\sim\operatorname{\as{ECR}(\beta,\lambda)} $ are defined by $ \sym{imS}(x_0)=\beta\lambda\int_0^{x_0}\nicefrac{x^{r+1}}{(\lambda^2+x^2)^{\nicefrac{3}{2}}}\left(1-\nicefrac{\lambda}{\sqrt{\lambda^2+x^2}}\right)^{\beta-1}\,\mathrm{d}x  $.
  Assuming the change of variable in \eqref{eq:ECR:cv1} and after some algebras we can express $ \sym{imS}(x_0)=\beta(\lambda\sqrt{2})^ru_0\,\int_0^{u_0}u^{\nicefrac{r}{2}+\beta-1}(1-u)^{-r}\left(1-\nicefrac{u}{2}\right)^{\nicefrac{r}{2}}\,\mathrm{d}u $,
  where $ u_0=1-\nicefrac{\lambda}{\sqrt{x_0^2+\lambda^2}} $. Now assuming $ t=\nicefrac{u}{u_0} $ we obtain
  \begin{equation}\label{eq:ECR:imsfinal}
    \sym{imS}(x_0)=\beta(\lambda\sqrt{2})^ru_0^{\nicefrac{r}{2}+\beta}\,\int_0^{1}t^{\nicefrac{r}{2}+\beta-1}(1-u_0t)^{-r}\left[1-\left(\frac{u_0}{2}\right)t\right]^{\nicefrac{r}{2}}\,\mathrm{d}t.
  \end{equation}
  Consider the following result
  \begin{equation}\label{eq:ECR:F1integral}
    \sym{F1}(a,b_1,b_2,c;x,y)=\frac{1}{\sym{bf}(a,c-a)}\,\int_0^1t^{a-1}(1-t)^{c-a-1}(1-xt)^{-b_1}(1-yt)^{-b_2}\,\mathrm{d}t,
  \end{equation}	
  unvielded by \citep[p.~77]{Bailey-Generalizedhypergeometricseries-1973}.
  The result is obtained by applying \cref{eq:ECR:F1integral} in \eqref{eq:ECR:imsfinal}.	
\end{proof}

\begin{proof}[Proof of \cref{prop:ECR:ros}]
  We need to determine the following integral
  \begingroup
  \allowdisplaybreaks
  \begin{align*}
    \sym{E}(X_{i:n}^r)=&\frac{\beta\lambda}{\sym{bf}(i,n-i+1)} \int_0^{\infty}x^r\frac{x}{(\lambda^2+x^2)^{\nicefrac{3}{2}}}\left(1-\frac{\lambda}{\sqrt{\lambda^2+x^2}}\right)^{i\beta-1}\nonumber\\
                       &\times\left[1-\left(1-\frac{\lambda}{\sqrt{\lambda^2+x^2}}\right)^{\beta}\right]^{n-i}\,\mathrm{d}x.
  \end{align*}
  \endgroup
  Using the change of variable in \eqref{eq:ECR:cv1} and after some algebras we can set
  \begin{equation}\label{eq:ECR:osstep1}
    \sym{E}(X_{i:n}^r)=\frac{\beta(\lambda\sqrt{2})^r}{\sym{bf}(i,n-i+1)}\int_0^1\frac{u^{\nicefrac{r}{2}+i\beta-1}(1-u)^{-r}}{(1-\nicefrac{u}{2})^{\nicefrac{-r}{2}}}(1-u^\beta)^{n-i}\,\mathrm{d}u.
  \end{equation}
  Note that we can express
  \begin{equation}\label{eq:ECR:binexp1}
    (1-u^\beta)^{n-i}=\sum_{j=0}^{n-i}(-1)^j\binom{n-i}{j}u^{j\beta}.
  \end{equation}
  Now applying \cref{eq:ECR:binexp1} in \eqref{eq:ECR:osstep1} we obtain
  \begin{equation*}
    \sym{E}(X_{i:n}^r)=\frac{\beta(\lambda\sqrt{2})^r}{\sym{bf}(i,n-i+1)}\sum_{j=0}^{n-i}(-1)^j\binom{n-i}{j}\int_0^1\frac{u^{\nicefrac{r}{2}+(i+j)\beta-1}(1-u)^{-r}}{(1-\nicefrac{u}{2})^{\nicefrac{-r}{2}}}\,\mathrm{d}u.
  \end{equation*}
  Finally using the identity in \eqref{eq:ECR:2F1integral} the result is obtained.
\end{proof}

\begin{proof}[Proof of \cref{prop:ECR:FEI}]
  Note that
  \begin{align}
    \sym{SC}_{\beta\beta}     & = -\frac{n}{\beta^2}\label{eq:ECR:FOIMele1}                                                                                                                                                                                    \\
    \sym{SC}_{\beta\lambda}   & = -\sum_{i=1}^{n}{\frac {\sqrt {{\lambda}^{2}+{x_i}^{2}}+\lambda}{{\lambda}^{2}+{x_i}^{2}} }\nonumber\label{eq:ECR:FOIMele2}                                                                                                   \\
    \sym{SC}_{\lambda\lambda} & = -\sum_{i=1}^{n}\frac{{x_i}^4+(\beta +4) \lambda ^2 {x_i}^2-\lambda ^3 \left[(\beta +1) \lambda +(\beta -1) \sqrt{\lambda ^2+{x_i}^2}\right]}{\lambda ^2 \left(\lambda ^2+{x_i}^2\right)^2}\nonumber
  \end{align}
  The result in \eqref{eq:ECR:FEIMele1} is derived taking the expectation of the constant obtained in \eqref{eq:ECR:FOIMele1}.
  The second element in \eqref{eq:ECR:FEIMele2} is obtained by the following integral.
  \begin{equation*}
    \sym{cum}_{\beta\lambda} = -n\beta\lambda\int_0^{\infty } \frac{x \left(\lambda +\sqrt{\lambda ^2+x^2}\right)}{\left(\lambda ^2+x^2\right) \left(\lambda ^2+x^2\right)^{\nicefrac{3}{2}}} \left(1-\frac{\lambda }{\sqrt{\lambda ^2+x^2}}\right)^{\beta -1}\,\mathrm{d}x.
  \end{equation*}
  Applying the transform in \eqref{eq:ECR:cv1} and then after some algebras we obtain	$ \sym{cum}_{\beta\lambda} = -\nicefrac{n\beta}{\lambda} \int_0^1 (u-2) (u-1) u^{\beta -1} \, \mathrm{d}u $,
  which is a simple integral whose value is given by $ \sym{cum}_{\beta\lambda} = -\nicefrac{n}{\lambda}\left[\frac{\beta +4}{\left(\beta+1\right)(\beta+2)}\right] $ and after partial fractions decomposition we obtain \eqref{eq:ECR:FEIMele2}.  
  The third element in \eqref{eq:ECR:FEIMele3} is provided by the following integral
  \begin{footnotesize}
    \begin{equation*}
      \sym{cum}_{\lambda\lambda} = -n\beta  \lambda  \int_0^{\infty } x \left\{\frac{ x^4+(\beta +4) \lambda ^2 x^2-\lambda ^3 \left[(\beta +1) \lambda +(\beta -1) \sqrt{\lambda ^2+x^2}\right] }{\lambda ^2 \left(\lambda ^2+x^2\right)^2 \left(\lambda ^2+x^2\right)^{\nicefrac{3}{2}}}\right\}\left(1-\frac{\lambda }{\sqrt{\lambda ^2+x^2}}\right)^{\beta -1}\,\mathrm{d}x.
    \end{equation*}
  \end{footnotesize}
  Applying the change of variable in \eqref{eq:ECR:cv1} and after some algebras we obtain
  \begin{equation*}
    \sym{cum}_{\lambda\lambda} = -\frac{n\beta }{\lambda ^2}\left\{\int_0^1 \left[2 (\beta +2) u^4-3 (3 \beta +5) u^3+(14 \beta +19) u^2-9 (\beta +1) u+2\beta\right] u^{\beta -1} \, \mathrm{d}u\right\},
  \end{equation*}
  which is a simple integral whose value is asserted by $ \sym{cum}_{\lambda\lambda} = -\nicefrac{n\beta  }{\lambda ^2}\left[\frac{\beta^2 +11\beta+36}{(\beta +2) (\beta +3) (\beta +4)}\right] $, and after partial fractions decomposition we obtain \eqref{eq:ECR:FEIMele3}.
\end{proof}

\begin{proof}[Proof of \cref{prop:ECR:TD}]
  Note that
  \begin{subequations}
    \begingroup
    \allowdisplaybreaks
    \begin{align}
      \sym{SC}_{\beta\beta\beta}       & = \frac{2n}{\beta ^3}\label{eq:ECR:TDbbb}                                                                                                                                               \\
      \sym{SC}_{\beta\beta\lambda}     & = 0\label{eq:ECR:TDbbl}                                                                                                                                                                 \\
      \sym{SC}_{\beta\lambda\lambda}   & = \sum_{i=1}^{n}\frac{\lambda  \left(\lambda +\sqrt{\lambda ^2+x_i^2}\right)-x_i^2}{\left(\lambda ^2+x_i^2\right)^2}\label{eq:ECR:TDbll}                                                \\
      \sym{SC}_{\lambda\lambda\lambda} & = \sum_{i=1}^{n}\frac{2 x_i^6+6 \lambda ^2 x_i^4-2 \lambda ^5 \left[(\beta +1) \lambda +(\beta -1) \sqrt{\lambda ^2+x_i^2}\right]}{\lambda ^3 \left(\lambda ^2+x_i^2\right)^3}\nonumber \\
                                       & \quad+\sum_{i=1}^{n}\frac{\lambda ^3 x_i^2 \left[6 (\beta +3) \lambda +(\beta -1) \sqrt{\lambda ^2+x_i^2}\right]}{\lambda ^3 \left(\lambda ^2+x_i^2\right)^3}\label{eq:ECR:TDlll}
    \end{align}
    \endgroup
  \end{subequations}
  The results in \eqref{eq:ECR:EVTDbbb} and \eqref{eq:ECR:EVTDbbl} are derived taking the expected values of the constants obtained in \eqref{eq:ECR:TDbbb} and \eqref{eq:ECR:TDbbl} respectively. The result in \eqref{eq:ECR:EVTDbll} is achieved by taking the expected value of \cref{eq:ECR:TDbll} which is determined by the integral
  \begin{equation*}
    \sym{cum}_{\beta\lambda\lambda}=-n\beta  \lambda  \int_0^{\infty } x\left[\frac{ \lambda  \left(\lambda +\sqrt{\lambda ^2+x^2}\right)-x^2}{\left(\lambda ^2+x^2\right)^2 \left(\lambda ^2+x^2\right)^{\nicefrac{3}{2}}}\right] \left(1-\frac{\lambda }{\sqrt{\lambda ^2+x^2}}\right)^{\beta -1}\,\mathrm{d}x.
  \end{equation*}
  Applying the change of variable in \eqref{eq:ECR:cv1} it is possible to rewrite this integral in the form $ \sym{cum}_{\beta\lambda\lambda}=\nicefrac{n\beta}{\lambda^2} \int_0^{1 } (u-2) (u-1)^2 (2 u-1) u^{\beta-1} \, \mathrm{d}u $, which integrand is polynomial and can be evaluated as $ \sym{cum}_{\beta\lambda\lambda}=-\nicefrac{2n}{\lambda^2}\left[\frac{ \beta^2 +4\beta-24}{(\beta +1) (\beta +2) (\beta +3) (\beta +4)}\right] $ and the result in \eqref{eq:ECR:EVTDbll} follows after partial fractions decomposition. The result in \eqref{eq:ECR:EVTDlll} is obtained by the following integral
  \begin{align*}
    \sym{cum}_{\lambda\lambda\lambda}&=n\beta  \lambda  \int_0^{\infty }x\left\{ \frac{2 x^6+6 \lambda ^2 x^4-2 \lambda ^5 \left[(\beta +1) \lambda +(\beta -1) \sqrt{\lambda ^2+x^2}\right]}{\lambda ^3 \left(\lambda ^2+x^2\right)^3 \left(\lambda ^2+x^2\right)^{\nicefrac{3}{2}}}\right\}\left(1-\frac{\lambda }{\sqrt{\lambda ^2+x^2}}\right)^{\beta -1} \\
                                     & \quad+  x\left\{ \frac{ \lambda ^3 x^2 \left[6 (\beta +3) \lambda +(\beta -1) \sqrt{\lambda ^2+x^2}\right]}{\lambda ^3 \left(\lambda ^2+x^2\right)^3 \left(\lambda ^2+x^2\right)^{\nicefrac{3}{2}}}\right\}\left(1-\frac{\lambda }{\sqrt{\lambda ^2+x^2}}\right)^{\beta -1}\,\mathrm{d}x
  \end{align*}
  Applying the change of variable in \eqref{eq:ECR:cv1} and after some algebras it is possible to rewrite this integral in the form
  \begin{align*}
    \sym{cum}_{\lambda\lambda\lambda}&=\frac{n\beta}{\lambda^3}\int_{0}^{1}\left[-8 (\beta +2) u^6+(51 \beta +93) u^5-3 (43 \beta +71) u^4+3 (55 \beta +81) u^3\right.\\
                                     &\quad\left.-3 (37 \beta +47) u^2+36 (\beta +1) u-4 \beta\right]u^{\beta-1}\,\mathrm{d}u,
  \end{align*}
  which integrand is polynomial and is given by $ \sym{cum}_{\lambda\lambda\lambda}=\nicefrac{2n}{\lambda^3} \left[\frac{ \beta  \left(\beta ^4+20 \beta ^3+158 \beta ^2+691 \beta +1866\right)}{(\beta +2) (\beta +3) (\beta +4) (\beta +5) (\beta +6)}\right] $ and the result in \eqref{eq:ECR:EVTDlll} follows after partial fractions decomposition.	
\end{proof}

\begin{proof}[Proof of \cref{prop:ECR:cox-snell-fisrtorderbias}]
  By \cref{eq:ECR:cox-snell-bias} we can express the second order bias of the \ab{MLE} as $ _{(\nicefrac{1}{n})}\sym{bias}(\hat{\theta}_i)=\sum_{r,s,t}\sym{cum}^{\theta_i,\theta_r}\sym{cum}^{\theta_s,\theta_t}\left(\sym{cum}_{\theta_r\theta_s}^{(\theta_t)}-\frac{1}{2}\sym{cum}_{\theta_r\theta_s\theta_t}\right)  $.
  The result is obtained applying \cref{coro:ECR:invFEI}, \cref{coro:ECR:FEIder} and \cref{prop:ECR:TD} in this equation. 
\end{proof}


\end{document}